\documentclass[a4paper,12pt]{article}
\usepackage[latin1]{inputenc}
\usepackage{amscd,amsmath,amssymb,amsfonts,latexsym,amsthm,mathrsfs}
\usepackage{fancyhdr}
\usepackage{float}
\usepackage{acronym}
\usepackage{graphicx,color,framed}
\usepackage{rcs}
\usepackage{caption}
\usepackage{slashbox}

\usepackage{natbib}

\input epsf.tex
\newcounter{example}[section]
\def\theexample{\arabic{example}}
{\noindent\\ \refstepcounter{example} {\it Example
\theexample\hspace{2pt}}
-\hspace{2pt} }{ 
 \begin{flushright}
 {\footnotesize $\square$}
 \end{flushright}
}

\title{Bayesian Model Selection for Beta Autoregressive Processes}
\author{Roberto Casarin\setcounter{footnote}{1}\footnotemark{}$\,$\hspace{5pt}
Luciana Dalla Valle\setcounter{footnote}{2}\footnotemark{}$\,$\hspace{5pt}
Fabrizio Leisen\setcounter{footnote}{3}\footnotemark{}$\,$\hspace{1pt}\setcounter{footnote}{4}\footnotemark{}\\\\
{\centering {\small \setcounter{footnote}{1}\footnotemark{} Università di Brescia}}\\
{\centering {\small \setcounter{footnote}{2}\footnotemark{} Università di Milano}}\\
{\centering {\small \setcounter{footnote}{3}\footnotemark{} Universidad Carlos III de Madrid}}}

\begin{document}

\maketitle

\renewcommand{\thefootnote}{\fnsymbol{footnote}}\setcounter{footnote}{1}
\footnotetext[5]{Corresponding author: {\tt fabrizio.leisen@gmail.com}. Corresponding address: Universidad Carlos III de Madrid,
Departamento de Estadística, Calle Madrid, 126, 28903 Getafe (Madrid), Spain.}

\renewcommand{\thefootnote}{\arabic{footnote}}\setcounter{footnote}{1}

\begin{abstract}
We deal with Bayesian inference for Beta autoregressive processes. We restrict our attention to the class of conditionally linear processes. These processes are particularly suitable for forecasting purposes, but are difficult to estimate due to the constraints on the parameter space. We provide a full Bayesian approach to the estimation and include the parameter restrictions in the inference problem by a suitable specification of the prior distributions. Moreover in a Bayesian framework parameter estimation and model choice can be solved simultaneously. In particular we suggest a Markov-Chain Monte Carlo (MCMC) procedure based on a Metropolis-Hastings within Gibbs algorithm and solve the model selection problem following a reversible jump MCMC approach.
\end{abstract}


\section{Introduction}
The analysis of time series data defined on a bounded interval (such as rates or proportions) has been a challenging issues for many years and still represents an open issue. For modelling data defined on a bounded interval there are at least two alternative approaches. Historically the main approach applies a transform to the data in order to map the interval to the real line and then uses standard time series models. Typical examples of transformations are the additive log-ratio transformation and the Box-Cox transformation (see \cite{Ait86}). One of the earlier and relevant contributions in this framework is \cite{Wal87}.

In this paper we follow the second approach, that is based on a direct modelling on the original sample space. Among the first contributions along this line we refer to \cite{GruRafGut93} who suggest a multivariate state space model for times series data defined on the standard simplex. Another seminal contribution is \cite{Mac85} which introduces a new Beta autoregressive process for times series defined on the standard unit interval $(0,1)$. In the recent years, \cite{FerCri04} introduce the general Beta regression model, showing that is more convenient to consider the data in the original sample space instead of using a transformation. In particular they use a reparametrization, which is often employed for the inference of Beta mixtures (see \cite{RobRou02}), that allows to identify the mean as a parameter of the distribution. \cite{RocCri09} extend the Beta regression model and propose a Beta autoregressive moving average process which possibly includes exogenous variables in the dynamics. In an applied context Beta regression models have been also recently used in \cite{AmiCas07} for modelling dynamic correlation, in \cite{Cal10} for modeling the recovery rate as a mixed random variable and in \cite{BilCas10a} and \cite{BilCas10b} for modeling the transition matrix of a Markov-switching model. In \cite{BilCas10b} a Bayesian procedure for latent Beta autoregressive models of the first order with nonlinear conditional mean is presented.

In particular, the main contribution of this paper is to propose a Bayesian estimation method that allows to estimate the parameters as well as the number of components of Beta autoregressive models. The parametrization used is the same given in \cite{RocCri09}. We handle the autoregressive part of the model, while the moving average part will be object of future research. Without loss of generality, we focus our attention on the linear conditional mean process, which is appealing for forecasting purposes but is more challenging due to the constraints on the parameters. The extension to the nonlinear conditional mean case requires minor modifications of the procedure proposed in this paper.
Considering the parameter estimation, we extend the work of \cite{BilCas10b} to the case of Beta autoregressive models of order $k$ and discuss the choice of priors for autoregressive coefficients, that becomes quite a challenging matter for a Beta autoregressive of general order.

Moreover, in order to select the model order, we propose a Reversible Jump Markov chain Monte Carlo (RJMCMC) algorithm for Beta autoregressive models, extending the works of \cite{BroGiuRob03} and \cite{EhlBro08} for Gaussian autoregressive models and \cite{EncNeaSub09} for Integer-Valued ARMA processes.

The reversible jump algorithm (see \cite{Gre95}) is an extension of the Metropolis-Hasting algorithm in which the dimension of the parameter space can vary between iterates of the Markov Chain. From a Bayesian point of view, suppose to have a set of models $\mathcal{M}=\lbrace\mathcal{M}_1,\mathcal{M}_2,\dots\rbrace$ where the model $\mathcal{M}_k$ has a vector of parameters {\bf $\theta_k$} of dimension $n_k$. The set $\mathcal{M}$ is indexed by a parameter $k\in\mathcal{K}$. Then the joint distribution of $(k,{\bf \theta_k})$ given the observed data {\bf $x$} is
\begin{equation}\label{Perve1}
\pi((k,{\bf \theta_k})\mid {\bf x})\propto L({\bf x}\mid (k,{\bf \theta_k}))p((k,{\bf \theta_k}))
\end{equation}
where $L({\bf x}\mid (k,{\bf \theta_k}))$ is the product of the likelihood and the joint prior $p((k,{\bf \theta_k}))=p({\bf \theta_k}|k)p(k)$ is constructed from the prior distribution of ${\bf \theta_k}$ under model $\mathcal{M}_k$ and the prior for the model
indicator $k$ (i.e. the prior for model $\mathcal{M}_k$).

The RJMCMC algorithm uses the joint posterior distribution in Equation (\ref{Perve1}) as the target of a Markov chain Monte Carlo sampler over the state space $\Theta=\cup_{k\in\mathcal{K}}(k,\mathbb{R}^{n_k})$ where the states of the Markov chain are of the form $(k,{\bf \theta_k})$, and the dimension of which can vary over the state space. Accordingly, from the output of a single Markov chain sampler, the user is able to obtain a full probabilistic description of the posterior probabilities of each model having observed the data, $x$, in addition to the posterior distributions of the individual models. For a recent account about reversible jump algorithm see Fan and Sisson (2010).

In this paper, in order to design an efficient reversible jump algorithm for Beta autoregressive processes we follow the strategy of \cite{BroGiuRob03} and \cite{EhlBro08}. We
consider jump move and stochastic reverse jump move between spaces of different dimensions and calibrate the parameters of the jump proposal distribution in order
to increase the acceptance rate of the move (see \cite{EhlBro08}). The Gaussian assumption allows \cite{EhlBro08} to have a close form solution of the calibration problem, while in the case of the Beta processes this result is no more valid. In our contribution we propose to combine an algorithm for the approximation of the parameter posterior mode with the calibration strategy of the jump proposal. The mode approximation allows us to find a close form solution to the parameter choice of the proposal distribution.

The outline of the paper is as follows. In Section 2 the Beta autoregressive process of order $k$ (BAR(k)) is introduced with a suitable parametrization. In Section 3 the Bayesian inference is developed for the BAR(k) model under the choice of some priors. In Section 4 the RJ algorithm is developed for the BAR processes, in particular
we studied two strategies. In Section 5 some simulation results are shown and in Section 6 the new machinery is applied to an unemployment rate and capacity utilization datasets.

\section{The Beta AR(k) model}
Let us define a Beta autoregressive process $\{x_{t}\}_{t\geq 0}$ of the
order $k$ as follows
\begin{equation}
x_{t}|\mathcal{F}^{X}_{t}\sim\mathcal{B}e(\eta_{1t},\eta_{2t})
\end{equation}
where the $\mathcal{F}^{X}_{t}=\sigma(\{x_{s}\}_{s\leq t})$ is the $\sigma$-algebra generated by the process,
$\mathcal{B}e(\eta_{1t},\eta_{2t})$ denotes the type I Beta distribution and $\eta_{1t}>0$ and $\eta_{2t}>0$
are the two parameters of the distribution, usually referred as shape parameters.
The two parameters are $\mathcal{F}_{t}$-measurable functions and depend possibly on the past values
of the process. The process has the following transition density
\begin{equation}
f(x_{t}|x_{t-1},\ldots,x_{t-k})=\frac{1}{B(\eta_{1t},\eta_{2t})}x_{t}^{\eta_{1t}-1}(1-x_{t})^{\eta_{2t}-1}\mathbb{I}_{(0,1)}(x_{t})
\end{equation}
where $B(a,b)$ is the Beta function with $a,b>0$ and $k$ is the order of the process. In the following we will denote
with $\hbox{BAR}(k)$ the $k$-order Beta autoregressive process and assume $k\leq k_{\max}$, with $k_{\max}<\infty$ the
maximum order of the process.

We considered the Beta distribution because it is a fairly flexible distribution. It is unimodal, uniantimodal,
increasing, decreasing or constant depending on the values of its parameters. We refer the reader to \citet{GupNad04}
for a review on the families of Beta distributions and to \citet{KotVan04}
for a review on the finite range distributions. We consider the Beta of the first type, because the first
and second moments of this distribution have a simple expression, which makes it easy to deal with the parameter
identification issues. In fact, as suggested by many parts in the literature, the parameters identification
should be accurately discussed. For inference and model interpretation purposes it is desiderable that the two
parameters $\eta_{1t}$ and $\eta_{2t}$ of the Beta distribution are not jointly determining both the shape and
the moments of the distribution.

In this work, in order to have an exact identification of the parameters of the conditional mean of the process and
of the shape parameter, we consider the parametrization of the Beta distribution suggested for example
in \citet{RobRou02} in the context of Bayesian inference for Beta mixtures and
in \citet{FerCri04} and \citet{RocCri09} in the context of
maximum likelihood inference for general beta regression models.

In the chosen parametrization the conditional distribution of the process at time $t$ is
\begin{equation}
x_{t}|\mathcal{F}_{t}\sim\mathcal{B}e(\eta_{t}\phi_{t},(1-\eta_{t})\phi_{t})
\end{equation}
where $\phi_{t}$ and $\eta_{t}$ represent the precision and location parameters respectively
(see \citet{RocCri09}). In this parametrization the conditional mean
of the process coincides with the location parameter and the variance is a quadratic function of
the location parameter as well
\begin{equation}
\mathbb{E}(x_{t}|\mathcal{F}_{t-1})=\eta_{t},\quad\mathbb{V}(x_{t}|\mathcal{F}_{t-1})=\eta_{t}(1-\eta_{t})/(1+\phi_{t})
\end{equation}
The process naturally exhibit heteroscedasticity, thus, in the following, we assume that the precision parameter is
constant over time, i.e. $\phi_{t}=\phi$.

We define $\boldsymbol{\alpha}=(\alpha_{0},\alpha_{1},\ldots,\alpha_{k})'$, $\mathbf{x}_{s:t}=(x_{s},\ldots,x_{t})'$, with $s<t$
and $\mathbf{z}_{t}=(1,\mathbf{x}_{t-1:t-k}')'$.
We assume that $\eta_{t}=\boldsymbol{\alpha}'\mathbf{z}_{t}$ is a linear combination of a constant and of
the $k$ past values of the process. We will assume that $\boldsymbol{\alpha}\in\Delta_{k+1}$ where $\Delta_{k+1}=\{\boldsymbol{\alpha}\in(0,1)^{k+1}|\sum_{i=0}^{k}\alpha_{i}\in(0,1)\}$ and refer to it as convexity constraints.

\begin{figure}[t]
\centering
\includegraphics[width=190pt]{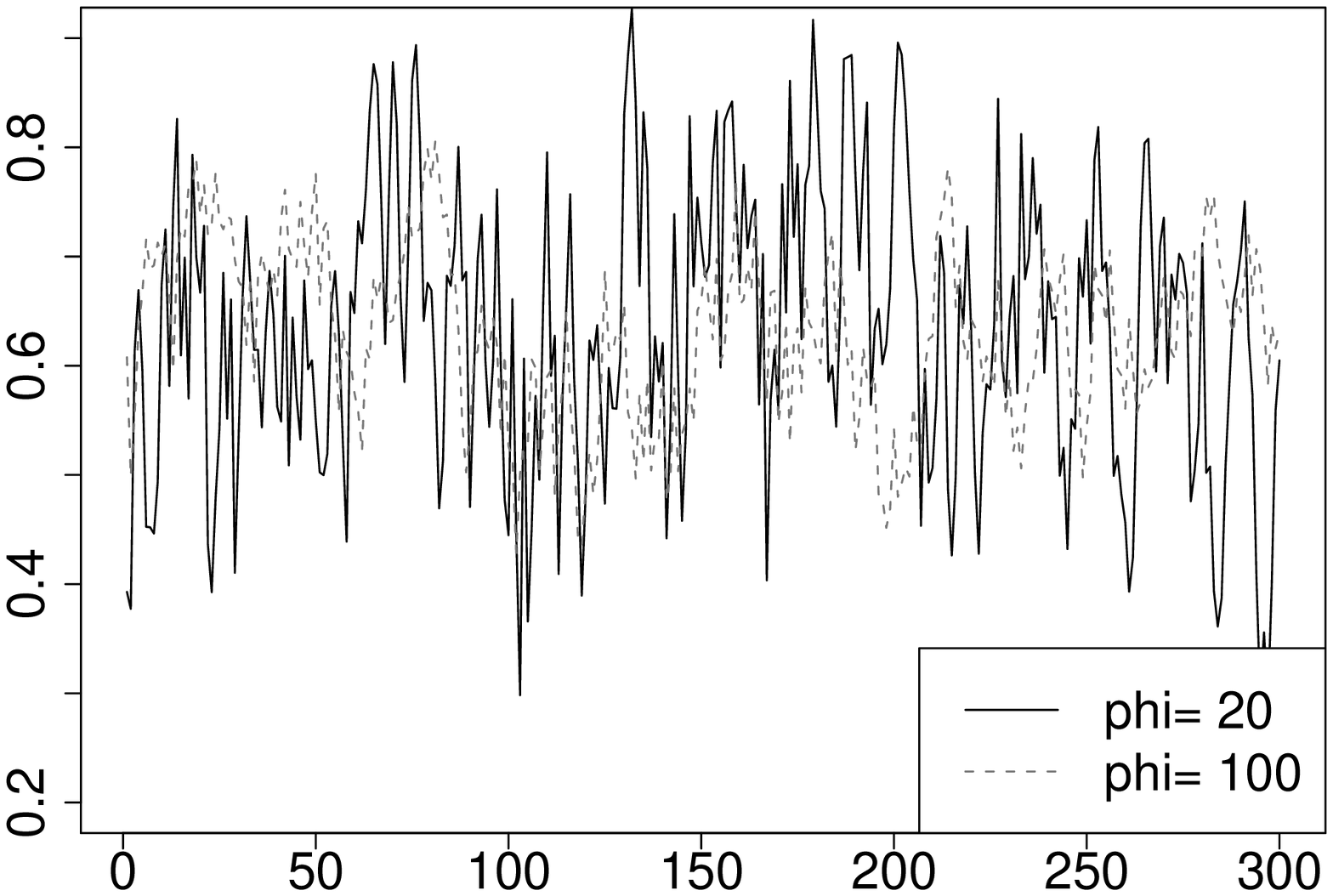}
\includegraphics[width=190pt]{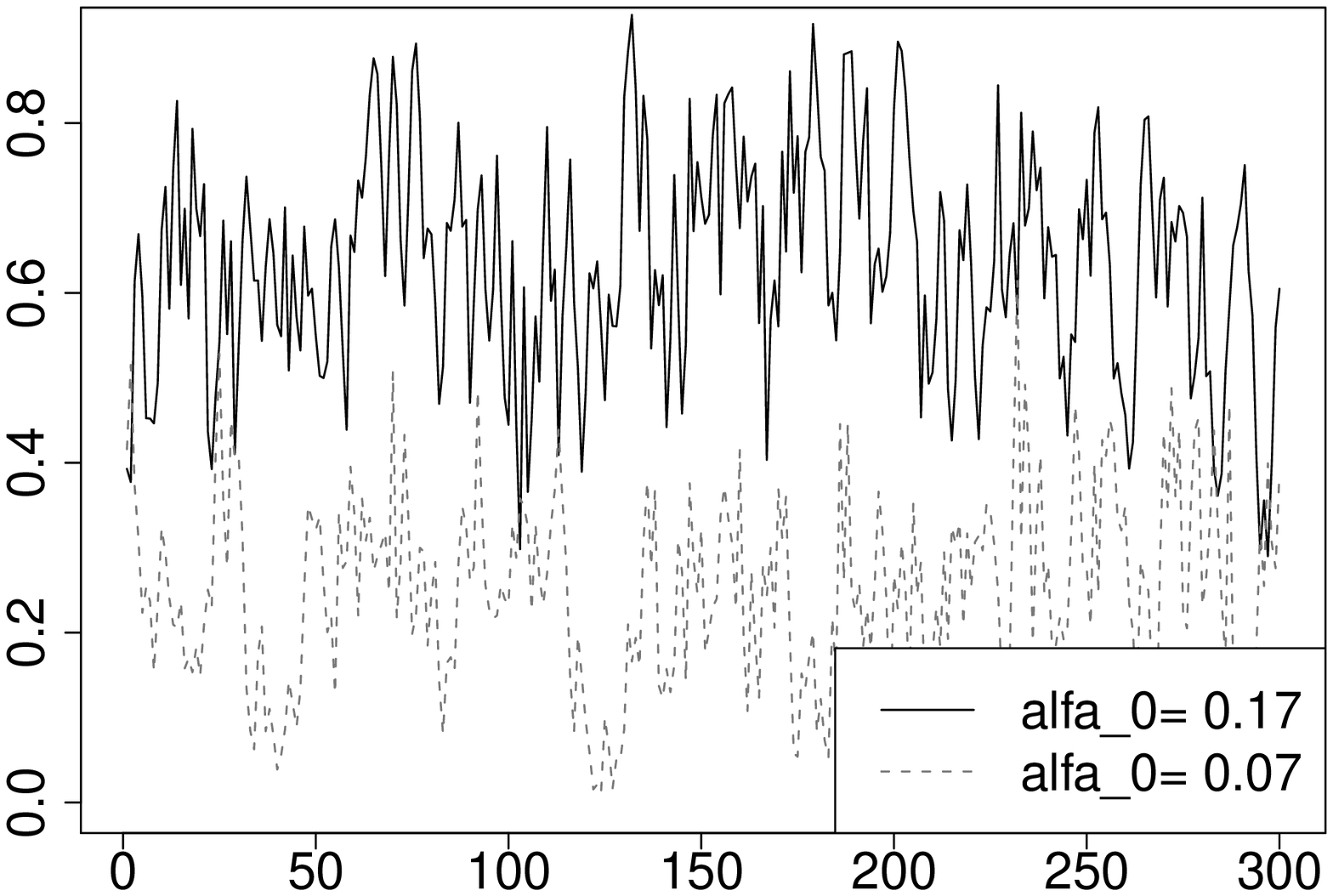}
\includegraphics[width=195pt]{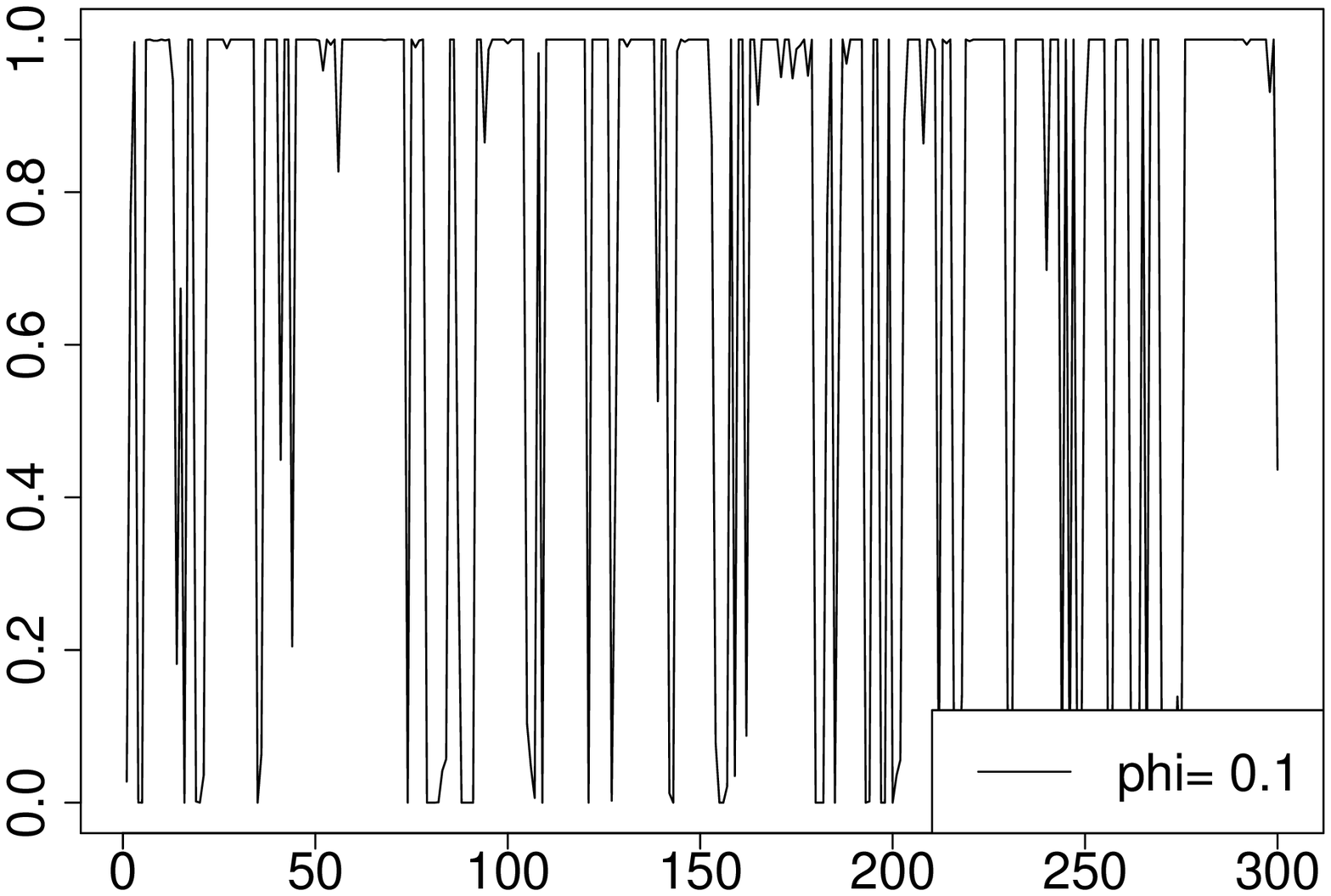}
\includegraphics[width=195pt]{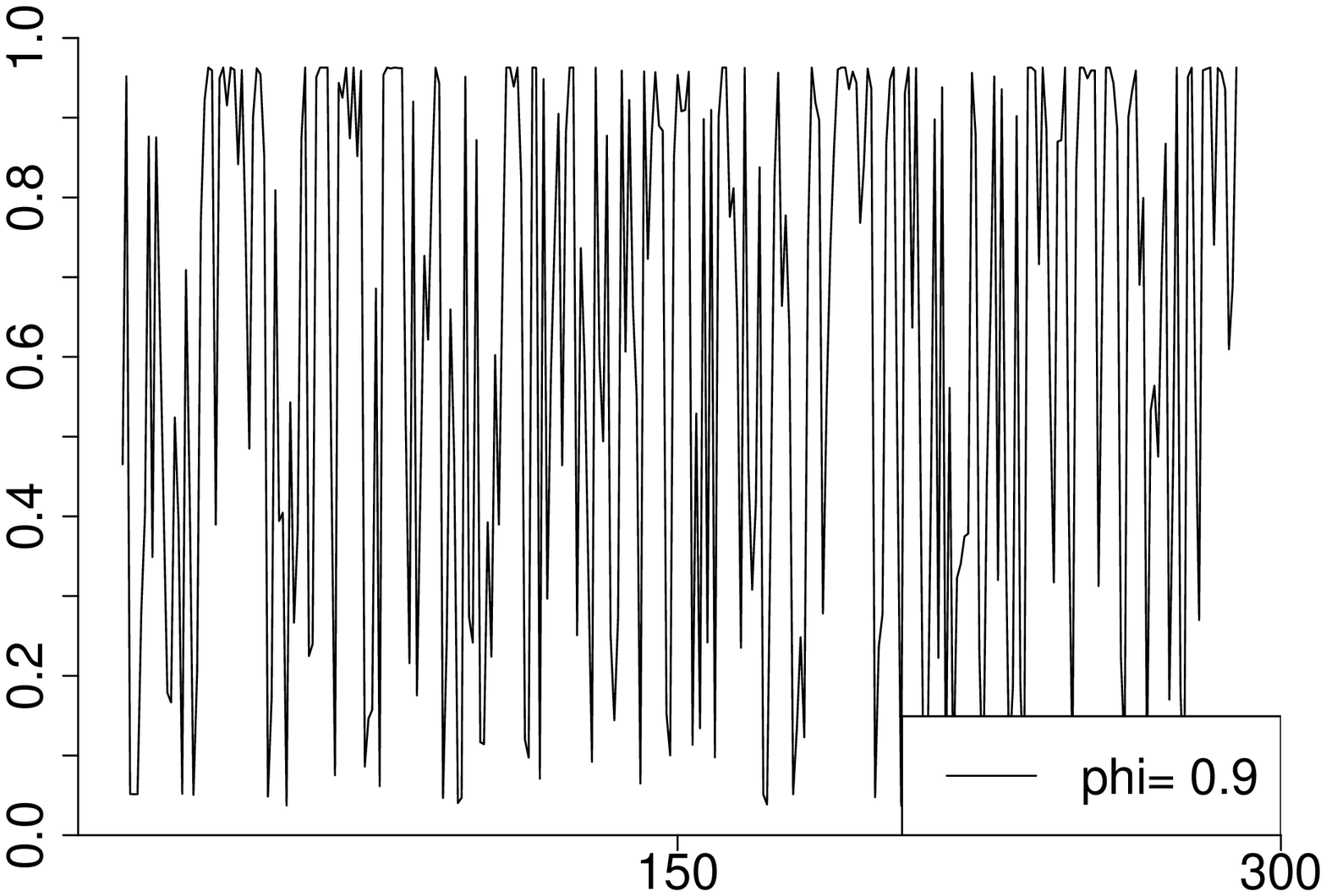}
\caption{Simulated trajectories of a $\hbox{BAR}(3)$ process for different parameter settings. \textit{Up-left}: the effect of the precision
parameter $\phi\in\{20,100\}$ for $\boldsymbol{\alpha}=(0.17,0.03, 0.1, 0.60)$. \textit{Up-right}: effect of the constant term $\alpha_{0}$ for
$\phi=100$ and $(\alpha_{1},\alpha_{2},\alpha_{3})=(0.03, 0.1, 0.60)$.
\textit{Bottom-left} and \textit{bottom-right}: anti-unimodal transition distribution and switching-type trajectories
of the BAR(3) process for different values of $\phi$ ($\phi\in\{0.1,0.9\}$) and with $\boldsymbol{\alpha}=(0.46,0.03, 0.01, 0.30)$.
Both of the cases correspond to $\underline{\eta}=0.46<1/\phi$ and $\bar{\eta}=0.8>(\phi-1)/\phi$.} \label{Fig_BARk}
\end{figure}

In Fig. \ref{Fig_BARk} there are some sample paths of a Beta process of the third order. The paths are given for different values of the
precision parameter $\phi$ and of the constant term $\alpha_{0}$. For higher values of $\phi$ the process exhibits
less volatility (upper-left chart). The larger the value of the constant term the greater is stationary mean of the
process (upper-right chart). Finally we observe that the conditional mean of the process $\eta_{t}\in[\underline{\eta},\bar{\eta}]$,
where $\underline{\eta}=\alpha_{0}$ and $\bar{\eta}=\alpha_{0}+\ldots+\alpha_{k}$. When $\underline{\eta}<1/\phi$ and
$\bar{\eta}>(\phi-1)/\phi$, then $\phi\eta_{t}<1$ and $\phi(1-\eta_{t})<1$ and the transition density of the process is anti-unimodal.
In this case the process exhibits a switching-type behavior (see the left and right charts at the bottom of Fig. \ref{Fig_BARk}).

This process extends to the order $k$ the first-order Beta process proposed in \citet{NieWal02}
and then discussed in a latent variable framework by \citet{AmiCas07}. The Beta process considered
in our work represents a special case of the beta regression model proposed in \citet{FerCri04}
and of the beta ARMA process introduced by \citet{RocCri09}. We focus on the linear case because it
is particularly suitable for forecasting purposes and is more difficult to estimate than the nonlinear case due to the
constraints on the parameters, which are necessary in order the process to be defined. In particular we show how to deal with the
first-order stationarity constraints in the inference procedure. The Bayesian inference approach presented in the following apply
in a straightforward way to the beta ARMA process with non-linear conditional mean given in \citet{RocCri09}.

\section{Bayesian inference}\label{Sec_Bayesian}
The likelihood function of the model is
\begin{equation}
\mathcal{L}(\boldsymbol{\alpha},\phi|\mathbf{x}_{t_{0}:T})=\prod_{t=t_{0}}^{T}
B(\eta_{t}\phi,(1-\eta_{t})\phi)^{-1}x_{t}^{\eta_{t}\phi-1}(1-x_{t})^{(1-\eta_{t})\phi-1}
\end{equation}
where $\mathbf{x}_{t_{0}:T}=(x_{t_{0}},\ldots,x_{T})'$ and $t_{0}=k_{\max}+1$. Note that we consider an
approximated likelihood because for a Beta process of the order $k\leq k_{\max}$
we assume the observations start in $t=k_{\max}+1$ and thus forget the
first $(k_{\max}-k)$ observations on $x_{t}$. Moreover in the following we will assume
that the first $k_{\max}$ initial values of the process are known. It is possible to include
the initial values in the inference process following, for example, the approach given
in \citet{VerAndDouGod04} for the Gaussian autoregressive processes.

\subsection{The priors}
The constant term and the coefficients of a $\hbox{BAR(k)}$ belong to the set $\Delta_{k+1}$, thus the usual assumption
of Gaussian prior distribution has to be modified. In the following we consider some alternative
prior specifications.

\begin{figure}[t]
\centering
\includegraphics[width=190pt, height=150pt,bbllx=20,bblly=215,bburx=567,bbury=620,clip=]{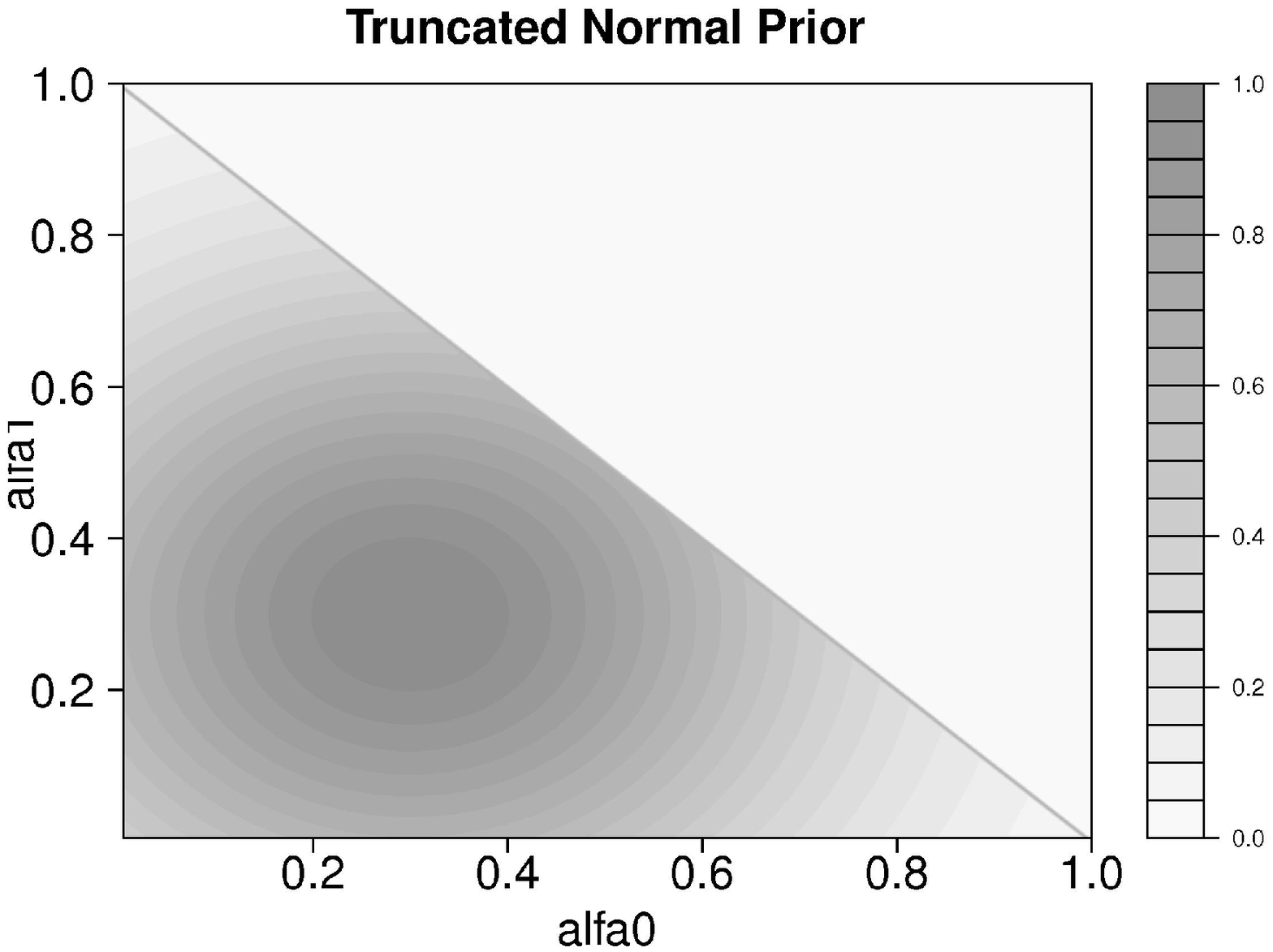}
\includegraphics[width=190pt, height=150pt,bbllx=20,bblly=215,bburx=567,bbury=620,clip=]{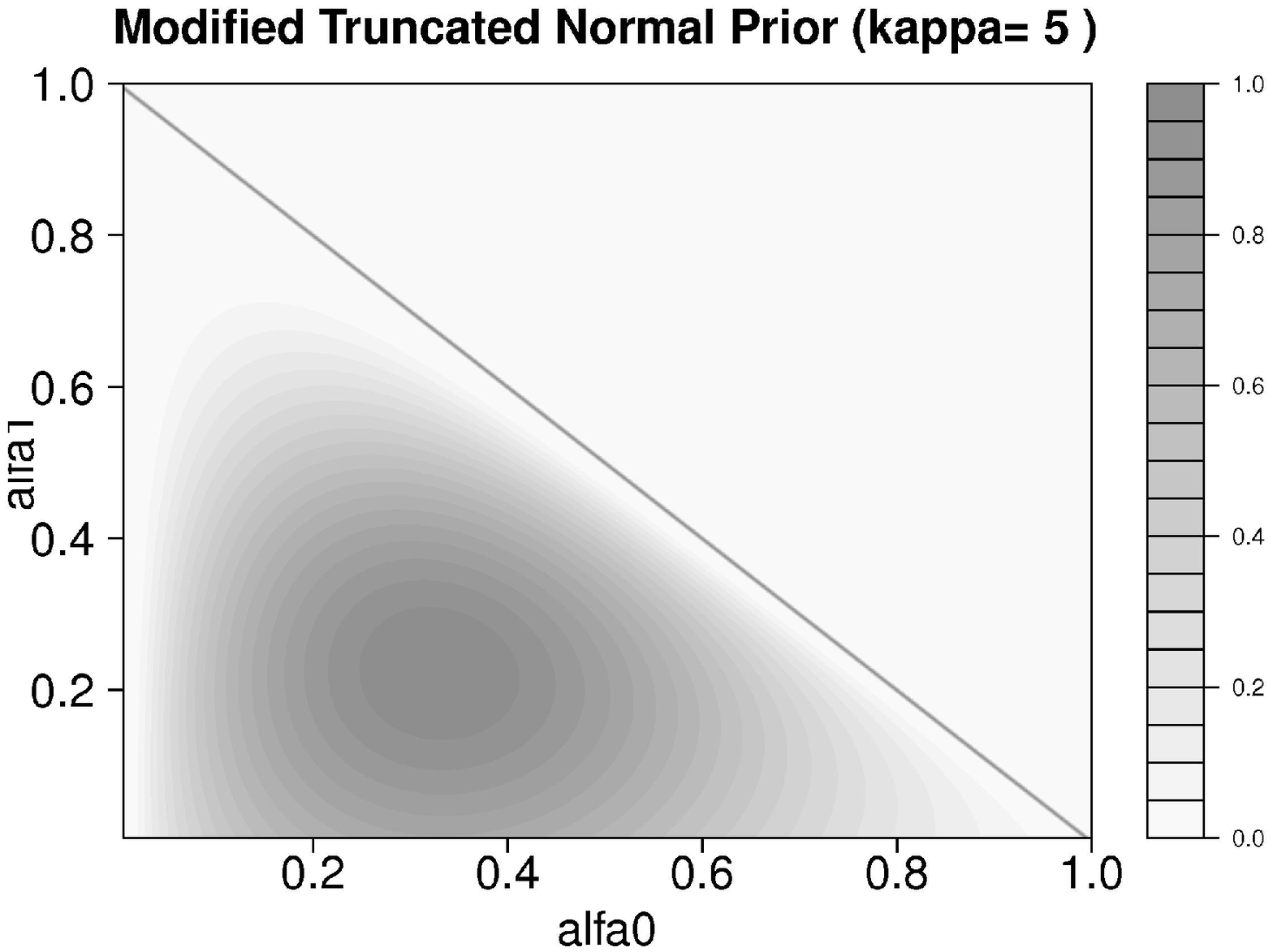}
\includegraphics[width=190pt, height=150pt,bbllx=20,bblly=215,bburx=567,bbury=620,clip=]{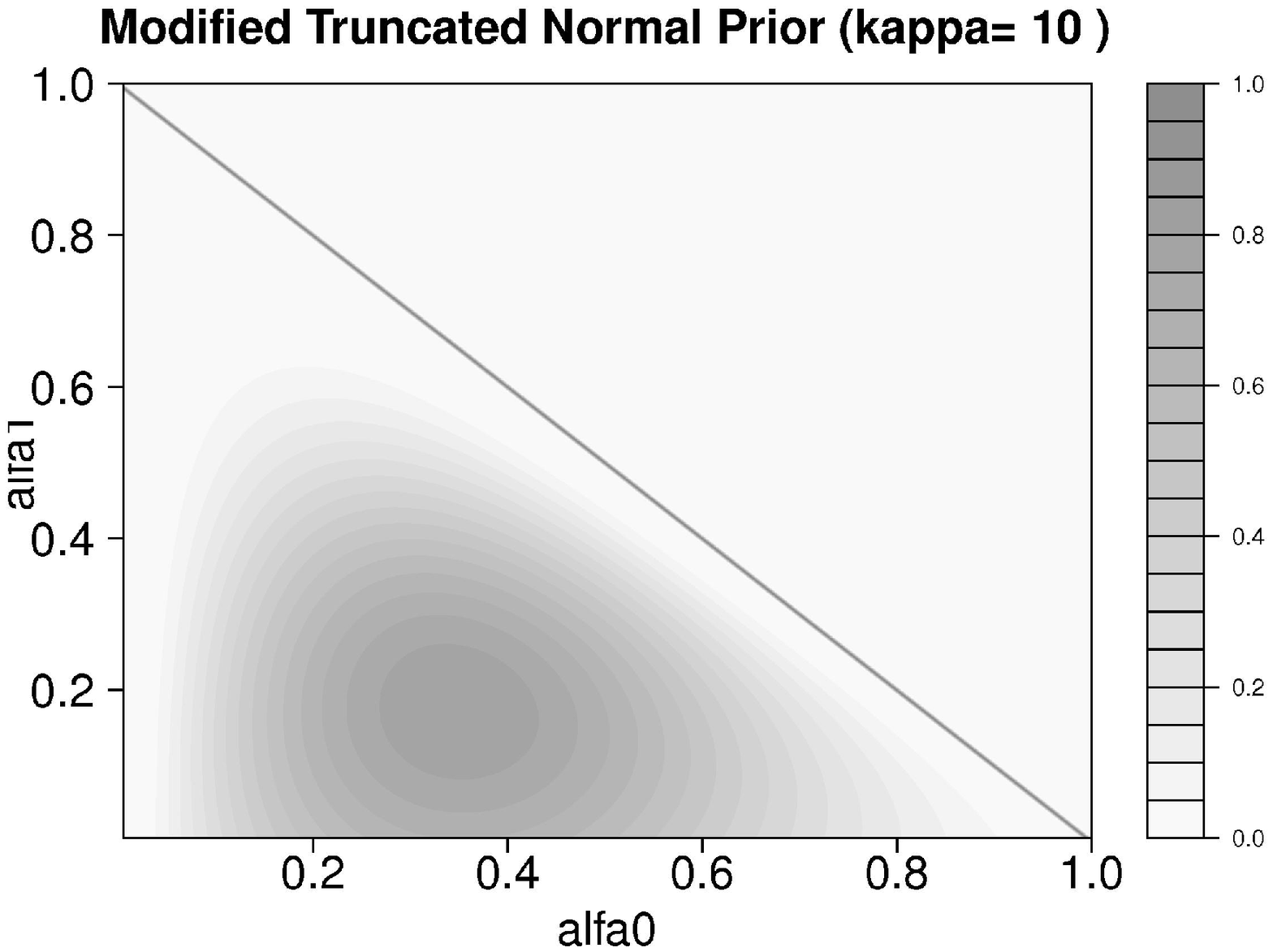}
\includegraphics[width=190pt, height=150pt,bbllx=20,bblly=215,bburx=567,bbury=620,clip=]{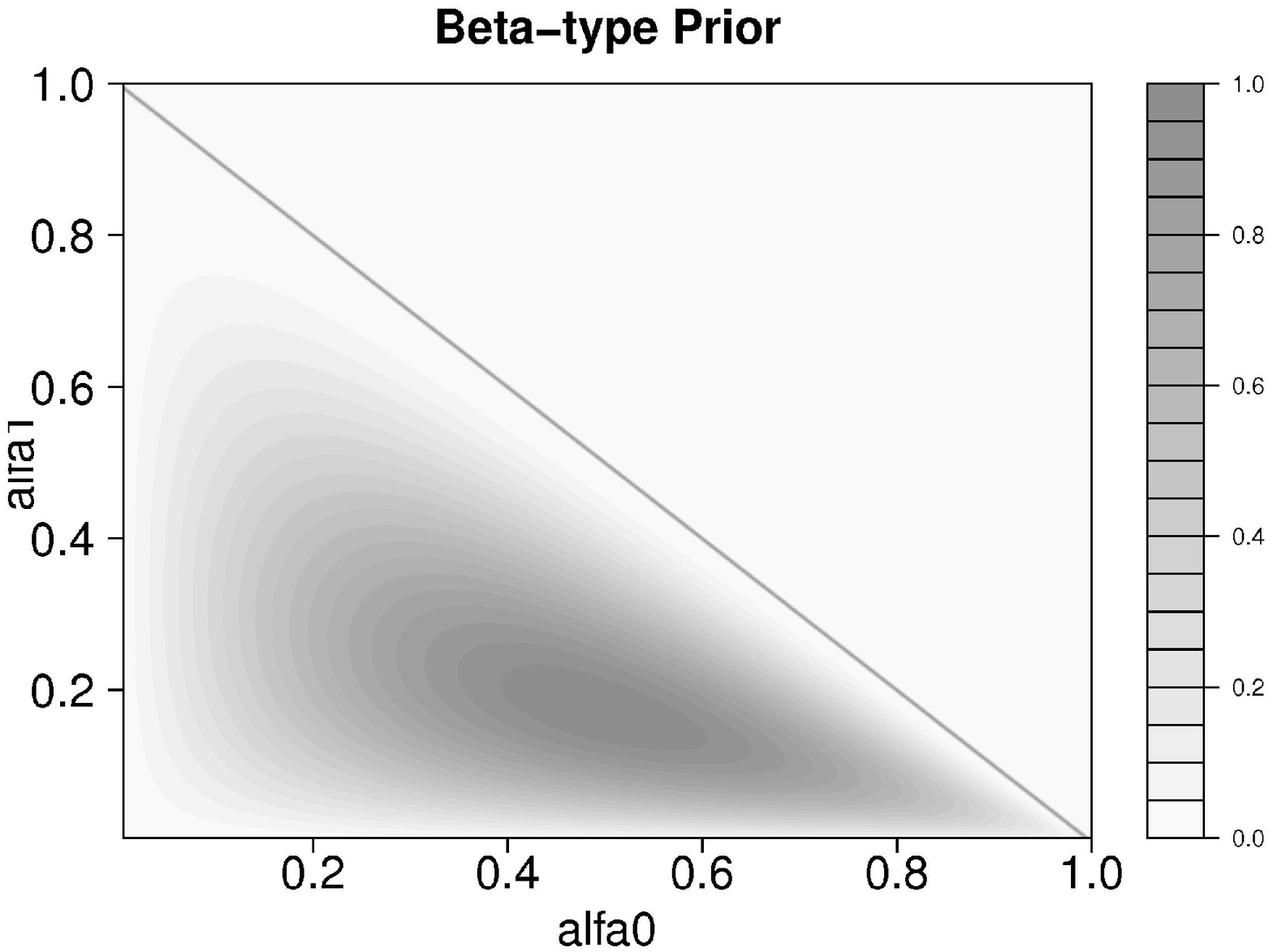}
\caption{Prior distributions for $\boldsymbol{\alpha}$ in the simplex $\Delta_{k+1}$ for $k=1$. In each graph,
the level set (\textit{gray areas}) and the sum-to-one constrain (\textit{solid line}). Truncated normal prior (\textit{up-left})
with parameters values $\boldsymbol{\nu}=(1/(k+2),1/(k+2))'$ and
$\Upsilon=0.1I_{2}$. Modified normal prior conditional on $\phi=10$, with parameters $\boldsymbol{\nu}=(1/(k+2),1/(k+2))'$
and $\Upsilon=0.1I_{2}$ for $\kappa=5$ (\textit{up-right}) and $\kappa=10$ (\textit{bottom-left}).
Multivariate Beta-type prior (\textit{bottom-right}) with parameters values $\boldsymbol{\nu}=(k+1,k+1)'$
and $\boldsymbol{\gamma}=(k+2,k+2)'$.} \label{Fig_Prior}
\end{figure}

We first considered a multivariate normal $\boldsymbol{\alpha}\sim \mathcal{N}_{k+1}\left(\boldsymbol{\nu},\Upsilon\right)$
with mean $\boldsymbol{\nu}$ and variance $\Upsilon$, truncated to the set $\Delta_{k+1}$. In the following we will denote with
\begin{equation}
f(\boldsymbol{\alpha})\propto\exp\left\{-\frac{1}{2}(\boldsymbol{\alpha}-\boldsymbol{\nu})'\Upsilon^{-1}(\boldsymbol{\alpha}-\boldsymbol{\nu})\right\}\mathbb{I}_{\Delta_{k+1}}(\boldsymbol{\alpha})
\end{equation}
the density function of the prior on $\boldsymbol{\alpha}$, where $\mathbb{I}_{A}(x)$ is the indicator function.
This prior distribution is given, for the case $k=1$, in the upper-left chart of Fig. \ref{Fig_Prior}. In the example
we assumed $\boldsymbol{\nu}=(1/(k+2),1/(k+2))'$ and $\Upsilon=0.1I_{2}$, where $I_{n}$ is the $n$-dimensional identity
matrix. Using this prior it is not easy to have a uniform prior and at the same time to guarantee the positivity of
the parameters of the Beta process transition density, i.e. $\eta_{t}\phi>0$ and $(1-\eta_{t})\phi>0$ $\forall t$.
The truncation on the simplex of the normal prior distribution can not guarantee that the inference procedure returns parameter values which satisfy at the desired constraints, and can generate numerical problems in the Monte Carlo procedures employed in the inference process.

In order to prevent the posterior to take values near the boundaries of the parameter space we follow \citet{RobRou02} and introduce a repulsive factor
around the boundaries of the standard simplex defined in the previous section.
We observe that $\underline{\eta}\leq\eta_{t}\leq\bar{\eta}$  and propose the following prior distribution conditional on $\phi$
\begin{equation}
f(\boldsymbol{\alpha})\propto
\exp\left\{-\frac{1}{2}(\boldsymbol{\alpha}-\boldsymbol{\nu})'\Upsilon^{-1}(\boldsymbol{\alpha}-\boldsymbol{\nu})\right\}
\exp\left\{-\frac{\kappa}{\phi^{2}\underline{\eta}(1-\bar{\eta})}\right\}
\mathbb{I}_{\Delta_{k+1}}(\boldsymbol{\alpha})
\end{equation}
where $\kappa$ is an hyperparameter. In the upper-right and bottom-left charts of Fig. \ref{Fig_Prior} we show, for the bivariate case, the shape of this prior distribution conditional on $\phi=10$ for two values of the hyperparameters ($\kappa=5$ and $\kappa=10$). The multiplicative
factor creates low density regions near the boundaries for all $t$. In both the simulation experiments and the real data applications
considered in Section \ref{Sec_SimRes}, this kind of prior contributes considerably to avoid the numerical problems in the
evaluation of the posterior density and of its gradient and Hessian which are needed for the Monte Carlo based inference procedures.
We note the solution of the numerical problems in the evaluation of the posterior distribution is due to the fact that the prior distribution
behaves like a penalizing term for the likelihood function and the penalization allows the algorithm to account for the constraints on the
parameters.
As an alternative, we consider a multivariate distribution which is naturally defined on $\Delta_{k+1}$ and which has density function
\begin{equation}
f(\boldsymbol{\alpha})\propto\prod_{i=0}^{k}\frac{1}{B(\nu_{i},\gamma_{i})}\left(\frac{\alpha_{i}}{A_{i}}\right)^{\gamma_{i}-1}\left(\frac{A_{i+1}}{A_{i}}\right)^{\nu_{i}-1}
\prod_{i=1}^{k}\frac{1}{A_{i}}\mathbb{I}_{\Delta_{k+1}}(\boldsymbol{\alpha})
\end{equation}
where $A_{i}=1-\sum_{j=0}^{i-1}\alpha_{i}$ and the parameters $\boldsymbol{\nu}=(\nu_{0},\ldots\nu_{k})'\in\mathbb{R}^{k+1}_{+}$
and $\boldsymbol{\gamma}=(\gamma_{0},\ldots,\gamma_{k})'\in\mathbb{R}^{k+1}_{+}$.

This distribution has been obtained by considering a set of independent Beta random variables
$v_{i}\sim{B}e(\nu_{i},\gamma_{i})$, with $i=0,\ldots,k$, and then the multivariate transformation
\begin{equation}
\quad \alpha_{j}=v_{j}\prod_{i=0}^{j-1}(1-v_{i}),\quad\hbox{for}\,\,j=0,\ldots,k
\end{equation}
This stochastic representation constitutes a natural way to generate random numbers from the multivariate
distribution defined above. It is easy to show that thanks to this representation the constraints
$\boldsymbol{\alpha}\in\Delta_{k+1}$ on the parameters of the Beta process are satisfied,
in particular $\sum_{i=0}^{k}\alpha_{i}=1-\prod_{i=0}^{k}(1-v_{i})\leq 1$. Bottom-right chart of
Fig. \ref{Fig_Prior} shows this prior for the parameter setting $\boldsymbol{\nu}=(k+1,k+1)'$
and $\boldsymbol{\gamma}=(k+2,k+2)'$. It can be seen from the picture that this prior gives a low probability mass to the values of $\boldsymbol{\alpha}$ near
the boundaries $\underline{\eta}=0$ and $\bar{\eta}=1$.

The precision parameter $\phi$ is positive, thus we assume a gamma prior
\begin{equation}
\phi\sim\mathcal{G}a(c,d)
\end{equation}
with parameters $c$ and $d$ and denote with $f(\phi)$ the associated density function.

\subsection{A MCMC algorithm}
In order to approximate the posterior mean for the parameters of the $\hbox{BAR}(k)$ process we consider a Gibbs algorithm. The full conditional distributions can not be simulated exactly and are simulated by a Metropolis-Hastings (M.-H.) step (see \citet{ChiGre95}). The full conditional distribution and the associated simulation method are described in the following.

The full conditional distribution of $\boldsymbol{\alpha}$ is
\begin{equation}
\pi(\boldsymbol{\alpha}|\phi,\mathbf{x}_{t_{0}:T})\propto\exp\left(
-\sum_{t=t_{0}}^{T}\log B(\eta_{t}\phi,(1-\eta_{t})\phi)+\sum_{t=t_{0}}^{T}A_{t}\eta_{t}\phi\right)f(\boldsymbol{\alpha})
\end{equation}
where $A_{t}=\log(x_{t}/(1-x_{t}))$. To simulate from this distribution we employ a M.-H. algorithm with a proposal distribution
which makes use of the information on the local structure of the posterior surface
(see \citet{AlbChi93} and \citet{LenDes00}). Consider the second-order Taylor expansion of the log-posterior, $g(\boldsymbol{\alpha})$,  centered around $\tilde{\boldsymbol{\alpha}}^{(j)}$.
\begin{equation}
g(\boldsymbol{\alpha})\approx g(\tilde{\boldsymbol{\alpha}}^{(j)})+
(\boldsymbol{\alpha}-\tilde{\boldsymbol{\alpha}}^{(j)})'\nabla^{(1)}g(\tilde{\boldsymbol{\alpha}}^{(j)})
+\frac{1}{2}(\boldsymbol{\alpha}-\tilde{\boldsymbol{\alpha}}^{(j)})'\nabla^{(2)}g(\tilde{\boldsymbol{\alpha}}^{(j)})(\boldsymbol{\alpha}-\tilde{\boldsymbol{\alpha}}^{(j)})
\end{equation}
where $\tilde{\boldsymbol{\alpha}}^{(j)}$ represents the approximated mode of the posterior, then at the $j$-th iteration of the
M.-H. step we generate a candidate as follows
\begin{equation}
\boldsymbol{\alpha}^{(*)}\sim\mathcal{N}(\tilde{\boldsymbol{\alpha}}^{(j)},\Sigma^{(j-1)})
\end{equation}
where $\Sigma^{(j-1)}=-\left(\nabla^{(2)}g(\tilde{\boldsymbol{\alpha}}^{(j-1)})\right)^{-1}$. The acceptance rate of this M.-H. step
can be easily found. Due to its relevance in the model selection procedure, the approximation method for the posterior mode
will be discussed in Section \ref{Sec_ModSec}.

The full conditional distribution of $\phi$ is
\begin{equation}
\!\pi(\phi|\boldsymbol{\alpha},\mathbf{x}_{t_{0}:T})\propto\exp\left(
\!-\!\!\sum_{t=t_{0}}^{T}\!\!\Big(\log B(\eta_{t}\phi,(1-\eta_{t})\phi)+\phi(A_{t}\eta_{t}+\log(1-x_{t}))\Big)\!\right)\!f(\phi)
\end{equation}
We simulate from the full conditional by a M.-H. step. We consider a Gamma random walk
proposal and at the $j$-th step of the algorithm, given the previous value
$\phi^{(j-1)}$ of the chain, we simulate
\begin{equation}
\phi^{(*)}\sim\mathcal{G}a(\sigma(\phi^{(j-1)})^{2},\sigma\phi^{(j-1)})
\end{equation}
and accept with probability
\begin{equation}
\min\left\{1,\frac{\pi(\phi^{(*)}|\boldsymbol{\alpha},\mathbf{x}_{t_{0}:T})}{\pi(\phi^{(j-1)}|\boldsymbol{\alpha},\mathbf{x}_{t_{0}:T})}
\frac{\Gamma(\sigma(\phi^{(j-1)})^{2})(\phi^{(j-1)})^{\sigma(\phi^{(*)})^2-1}(\sigma\phi^{(*)})^{\sigma(\phi^{(*)})^{2}}}
{\Gamma(\sigma(\phi^{(*)})^{2})(\phi^{(*)})^{\sigma(\phi^{(j-1)})^2-1}(\sigma\phi^{(j-1)})^{\sigma(\phi^{(j-1)})^{2}}}\right\}
\end{equation}

\section{Model Selection}\label{Sec_ModSec}
In this work we propose a Reversible-Jump MCMC (RJMCMC) approach to model selection for Beta autoregressive.
See \citet{BroGiuRob03} for an introduction to efficient proposal design
in RJMCMC algorithms and \citet{EhlBro08} for an application of some efficient
adaptive proposals to Gaussian autoregressive processes. See also \citet{VerAndDouGod04}
for an alternative RJMCMC schemes for Gaussian autoregressive. To the best of our knowledge only a few studies exist on the application of RJMCMC algorithms to non-Gaussian autoregressive models. Among the other we refer the interested
reader to \citet{EncNeaSub09} who proposed a RJMCMC algorithm for model selection
in integer valued ARMA models.

In this work we extend the RJMCMC approach of \citet{EhlBro08} to
the case of the Beta autoregressive processes of a unknown order. We consider the definition
of Beta autoregressive given in \citet{RocCri09} and restrict our attention
to the case of conditionally linear processes because the inference in this context is more difficult
due to the restriction on the parameters. Nevertheless the proposed RJMCMC algorithm can be easily
extended to the case of a Beta process with a nonlinear conditional mean.

In the following we will show how to deal, in a Beta autoregressive context, with some special constraints
which can be imposed in the inference process.  First we show how to deal with first order stationarity constraints. We suggest
to move between $\hbox{BAR}(k)$ and $\hbox{BAR}(k+1)$ processes and to consider the
innovation process naturally associated to a $\hbox{BAR}(k)$ to design this move. As suggested in \citet{EhlBro08}
we choose the parameters of the proposal for jumps between $\hbox{BAR}(k+1)$ and $\hbox{BAR}(k)$ in
such a way to make the acceptance rate close to one.

Secondly we consider convexity constraints on the parameters of the conditional means and suggest to use a trans-dimensional MCMC chain
which moves between spaces of any finite dimensions $k$ and $k'$ by employing a suitable proposal distribution.
In order to make the acceptance rate close to one we suggest to combine the second-order proposal strategy with a Newton-Rapson type sequential approximation of the posterior mode.

\subsection{Stationarity Constraints}
Let us define the innovation process $\xi_{t}=x_{t}-\eta_{t}$, then
\begin{equation}
\xi_{t}=(1-\alpha_{1}L-\ldots-\alpha_{k}L^{k})x_{t}-\alpha_{0}
\end{equation}
where $L$ is the lag operator. The innovation process can be written in terms of reciprocal roots
\begin{equation}
\xi_{t}=\prod_{j=1}^{k}(1-\lambda_{j}L)x_{t}-\alpha_{0}
\end{equation}
and it is stationary in mean if $|\lambda_{j}|<1$.

Let $\xi_{t}$ be the innovation term associated to $\hbox{BAR}(k)$ and $\xi_{t}^{*}$ the one associated to
a $\hbox{BAR}(k+1)$ then
\begin{eqnarray}
\xi_{t}^{*}&=&(1-rL)\prod_{j=1}^{k}(1-\lambda_{j}L)x_{t}-\alpha_{0}^{*}\nonumber\\
&=&(1-rL)(\xi_{t}+\alpha_{0})-\alpha_{0}^{*}\nonumber\\
&=&(1-rL)\xi_{t}
\end{eqnarray}
if $\alpha_{0}=\alpha_{0}^{*}$.

Thanks to this representation of the location parameter, the likelihood
function of the $\hbox{BAR}(k+1)$ process will write as
\begin{eqnarray}
&&\mathcal{L}(\alpha_{0},\boldsymbol{\lambda},r,\phi|\mathbf{x}_{t_{0}:T})=\prod_{t=t_{0}}^{T}
B(\eta_{t}^{*}\phi,(1-\eta_{t}^{*})\phi)^{-1}x_{t}^{\eta_{t}^{*}\phi-1}(1-x_{t})^{(1-\eta_{t}^{*})\phi-1}\nonumber\\
&&=\prod_{t=t_{0}}^{T}
B(g_{t}(r)\phi,(1-g_{t}(r))\phi)^{-1}
x_{t}^{g_{t}(r)\phi-1}(1-x_{t})^{(1-g_{t}(r))\phi-1}
\end{eqnarray}
where $g_{t}(r)=x_{t}-(\xi_{t}-r\xi_{t-1})$,
$\boldsymbol{\lambda}=(\lambda_{1},\ldots,\lambda_{k})'$ is the vector of
reciprocal roots of the $\hbox{BAR}(k)$ process which coincide with the first $k$ reciprocal
roots of the $\hbox{BAR}(k+1)$ and $r\in(-1,1)$ is the $k+1$-th reciprocal root
of the $\hbox{BAR}(k+1)$ process.

We assume independent Beta priors for the reciprocal roots of the process and in order to assure that the constraints on the
autoregressive coefficients are satisfied (see Appendix A), we let $u_{j}\sim\mathcal{B}e(a,b)$ and
define $\lambda_{j}=2u_{j}-1$ for $j=1,\ldots,k+1$, with $\lambda_{k+1}=r$. We denote with $f(\boldsymbol{\lambda})$
and $f(r)$ the densities associated to these priors.

In order to design a reversible jump algorithm, which allows for jumps between posterior distributions with
different dimensions, we need a model prior, a proposal distribution for the model order and for the new
parameters and finally a link function between the parameter
spaces of the two models. In the following we will assume that the model prior associated to $k$ is uniform:
$k\sim\mathcal{U}_{\{1,\ldots,k_{\max}\}}$, where $k_{\max}$ is the maximum order of the process and denote with $f(k)$
the associated density function. As a proposal $p_{k,k+1}$ for the order of the new model we consider a Bernoulli
distribution with probability $1/2$. We assume a truncated normal distribution as a proposal for the new root, i.e. $r\sim\mathcal{N}(\mu,\sigma^{2})$ truncated to $(-1,1)$ and denote with $q(r)$ its density.
Finally, as suggested in \citet{EhlBro08}, we consider an identity dimension matching function $f(\boldsymbol{\lambda},\lambda_{k+1})=(\boldsymbol{\lambda}',r)'$.

The acceptance ratio of the move from a $k$-order autoregressive to an autoregressive of order $(k+1)$ is
\begin{equation}
A_{k,k+1}=\frac{\mathcal{L}(\alpha_{0},\boldsymbol{\lambda},r,\phi|\mathbf{x}_{t_{0}:T})f(\alpha_{0})f(\boldsymbol{\lambda})f(\phi)f(r)f(k+1)p_{k+1,k}}
{\mathcal{L}(\alpha_{0},\boldsymbol{\lambda},\phi|\mathbf{x}_{t_{0}:T})f(\alpha_{0})f(\boldsymbol{\lambda})f(\phi)f(k)p_{k,k+1}}
\frac{1}{q(r)}|J|
\end{equation}
where $|J|$ is the Jacobian of the transformation. Under the above assumptions the acceptance rate simplifies to
\begin{equation}
A_{k,k+1}=\frac{\mathcal{L}(\alpha_{0},\boldsymbol{\lambda},r,\phi|\mathbf{x}_{t_{0}:T})f(r)}
{\mathcal{L}(\alpha_{0},\boldsymbol{\lambda},\phi|\mathbf{x}_{t_{0}:T})}
\frac{1}{q(r)}
\end{equation}

As suggested in \citet{EhlBro08} the likelihood of the
$\hbox{BAR}(k)$ and the $\hbox{BAR}(k+1)$ models are identical at $r=0$, which
is a weak non-identifiability centring point. We use the centring point to find the parameters of the
second-order proposal. To this aim we consider the first and second order derivatives of the
log-acceptance ratio.

The log-acceptance ratio of the move from a Beta process of order $k$ to one of order $k+1$ is
\begin{eqnarray}
&&\log A_{k,k+1}=C+\frac{1}{2\sigma^{2}}(r-\mu)^{2}-\sum_{t}\left(\log B(g_{t}(r)\phi,(1-g_{t}(r))\phi)\right)\nonumber\\
&&\quad\quad+\sum_{t=t_{0}}^{T}\left((g_{t}(r)\phi-1)\log(x_{t})\right)+\sum_{t}\left(((1-g_{t}(r))\phi-1)\log(1-x_{t})\right)\nonumber\\
&&\quad\quad+(a-1)\log(1+r)+(b-1)\log(1-r)
\end{eqnarray}
where $C$ is the log of a normalizing constant which does not depend on $r$.
The first- and second-order derivatives of the log-acceptance ratio are
\begin{eqnarray}
&&\partial_{r}\log A_{k,k+1}=\sum_{t=t_{0}}^{T}\xi_{t-1}\phi\left(\Psi^{(0)}((1-g_{t}(r))\phi)-\Psi^{(0)}(g_{t}(r)\phi)\right)\nonumber\\
&&\quad\quad+\sum_{t=t_{0}}^{T}\xi_{t-1}\phi\left(\log(x_{t})-\log(1-x_{t})\right)
+\frac{r-\mu}{\sigma^{2}}+\frac{a-1}{1+r}-\frac{b-1}{1-r}\\
&&\partial_{rr}\log A_{k,k+1}=-\sum_{t=t_{0}}^{T}\xi_{t-1}^{2}\phi^{2}\left(\Psi^{(1)}(g_{t}(r)\phi)+
\Psi^{(1)}((1-g_{t}(r))\phi)\right)\nonumber\\
&&\quad\quad+\frac{1}{\sigma^{2}}-\frac{a-1}{(1+r)^{2}}-\frac{b-1}{(1-r)^{2}}
\end{eqnarray}
where $\Psi^{(0)}$ and $\Psi^{(1)}$ are the digamma and trigamma functions respectively.

We evaluate the first- and second-order derivatives at $r=0$, solve for $\mu$ and $\sigma^{2}$
and find the parameters of the proposal
\begin{eqnarray}
\mu&=&\frac{\left((a-b)+\phi(U_{2}+U_{1})\right)}{(a+b)-2+\phi^{2}U_{3}}\\
\sigma^{2}&=&\frac{1}{(a+b)-2+\phi^{2}U_{3}}
\end{eqnarray}
where
\begin{eqnarray*}
U_{1}&=&\sum_{t=t_{0}}^{T}\xi_{t-1}\left(\Psi^{(0)}((1-x_{t}+\xi_{t})\phi)-\Psi^{(0)}((x_{t}-\xi_{t})\phi)\right)\\
U_{2}&=&\sum_{t=t_{0}}^{T}\xi_{t-1}\log(x_{t}/(1-x_{t}))\\
U_{3}&=&\sum_{t=t_{0}}^{T}\xi_{t-1}^{2}\left(\Psi^{(1)}((x_{t}-\xi_{t})\phi)+
\Psi^{(1)}((1-x_{t}+\xi_{t})\phi)\right)
\end{eqnarray*}

\subsection{Convexity Constraints}
In order to deal with the convexity constraints on the parameters in the conditional mean of the BAR process we suggest to use RJMCMC moves between spaces with any positive dimensions difference, and to control for the constraints by a suitable choice of the proposal distribution for the parameters.

Given the parameter vector $\boldsymbol{\alpha}\in\Delta_{k+1}$, the RJMCMC proposes a jump to the space $\Delta_{k'+1}$ by generating $k'$ from a discrete distribution $p_{k,k'}$ and then proposing a $k'$-dimensional parameter vector $\mathbf{u}\in\Delta_{k'+1}$ from the distribution $q(\mathbf{u})$, which we assumed to be defined on the standard simplex. For the dimensions matching we consider the transform $(\boldsymbol{\alpha}',\mathbf{u}')=(\mathbf{u},\boldsymbol{\alpha})$, which has unitary Jacobian. The acceptance probability of the proposal is
\begin{equation*}
A_{k,k'}=\frac{\mathcal{L}(\mathbf{u},\phi|\mathbf{x}_{t_{0}:T})f(\mathbf{u})f(\phi)p_{k'}p_{k',k}}
{\mathcal{L}(\boldsymbol{\alpha},\phi|\mathbf{x}_{t_{0}:T})f(\boldsymbol{\alpha})f(\phi)p_{k'}p_{k,k'}}\frac{q(\boldsymbol{\alpha})}{q(\mathbf{u})}
\end{equation*}
As proposal we consider a parametric family of distributions and choose the parameters in such a way that the log-acceptance rate is approximately equal to zero.

We consider a second-order methods and focus on the gradient and the Hessian of the log-acceptance rate with respect to $(\mathbf{u},\boldsymbol{\alpha})$. Note that the gradient naturally splits in the two subvectors, that are the gradient with respect to $\mathbf{u}$ and the one with respect to $\boldsymbol{\alpha}$. The cross derivatives in the Hessian are null. In the following we consider the case for $\mathbf{u}$, the case for $\boldsymbol{\alpha}$ being similar.

We choose the parameters of the proposal such that
\begin{eqnarray*}
\nabla^{(1)}\log A_{k,k'}&=&\sum_{t=t_{0}}^{T}\left(A_{t}-\Psi^{(0)}(\eta_{t}\phi)+\Psi^{(0)}((1-\eta_{t})\phi)\right)\phi\mathbf{z}_{t}\\
&\quad&+\nabla^{(1)}\log f(\mathbf{u})-\nabla^{(1)}\log q(\mathbf{u})=\mathbf{0}\\
\nabla^{(2)}\log A_{k,k'}&=&-\sum_{t=t_{0}}^{T}\left(\Psi^{(1)}(\eta_{t}\phi)+\Psi^{(1)}((1-\eta_{t})\phi)\right)\phi^{2}\mathbf{z}_{t}\mathbf{z}_{t}'\\
&\quad&+\nabla^{(2)}\log f(\mathbf{u})-\nabla^{(2)}\log q(\mathbf{u})=0
\end{eqnarray*}
where $\Psi^{(0)}$ and $\Psi^{(1)}$ have been defined in the previous section.

It should be noticed that as opposite to the Gaussian autoregressive model studied in \citet{EhlBro08}, the gradient and the Hessian depend on $\mathbf{u}$. In this work we decide to evaluate the derivatives at the approximated posterior mode $\tilde{\mathbf{u}}$, which is defined on the $k'$-dim space.

As an example let use consider a Gaussian distribution as proposal with mean $\boldsymbol{\mu}_{k'}$ and variance $\Sigma_{k'}$. Assume that the prior is the Gaussian distribution truncated on the simplex, then we solve the above system of equations (see Appendix B for the gradient and the Hessian in the Gaussian case) with respect to $\boldsymbol{\mu}_{k'}$ and $\Sigma_{k'}$ and find
\begin{eqnarray*}
\Sigma^{-1}_{k'}&=&\Upsilon_{k'}^{-1}+\sum_{t=t_{0}}^{T}\left(\Psi^{(1)}(\eta_{t}\phi)+\Psi^{(1)}((1-\eta_{t})\phi)\right)\phi^{2}\mathbf{z}_{t}\mathbf{z}_{t}'\\
\boldsymbol{\mu}_{k'}&=&\tilde{\mathbf{u}}+\Sigma_{k'}\left(\Upsilon_{k'}^{-1}(\tilde{\mathbf{u}}-\boldsymbol{\nu}_{k'})-\sum_{t=t_{0}}^{T}\left(A_{t}-\Psi^{(0)}(\eta_{t}\phi)+\Psi^{(0)}((1-\eta_{t})\phi)\right)\phi\mathbf{z}_{t}\right)
\end{eqnarray*}
By applying the same procedure for the derivatives with respect to $\boldsymbol{\alpha}$ we obtain a full modal
centering procedure, which promotes jumps between modes in the $k$-dimensional and $k'$-dimensional spaces.

We do not know the value of the posterior mode and in order to find an approximation we suggest to use a Newton-Rapson
procedure. More specifically at the iteration $j$ of the RJMCMC procedure the approximated mode $\tilde{\mathbf{u}}^{(j)}$ is determined by the following recursion
\begin{equation}
\tilde{\mathbf{u}}^{(j)}=\tilde{\mathbf{u}}^{(j-1)}-\Sigma^{(j-1)}\nabla^{(1)}g(\tilde{\mathbf{u}}^{(j-1)})
\end{equation}
where $\Sigma^{(j-1)}=-\left(\nabla^{(2)}g(\tilde{\mathbf{u}}^{(j-1)})\right)^{-1}$. The gradient and the Hessian of the log-posterior, $\nabla^{(1)}g$ and $\nabla^{(2)}g$ respectively, are
\begin{eqnarray*}
&&\nabla^{(1)}g(\mathbf{u})=\sum_{t=t_{0}}^{T}\left(A_{t}-\Psi^{(0)}(\eta_{t}\phi)+\Psi^{(0)}((1-\eta_{t})\phi)\right)\phi\mathbf{z}_{t}+\nabla^{(1)}\log f(\mathbf{u})\\ 
&&\nabla^{(2)}g(\mathbf{u})=-\sum_{t=t_{0}}^{T}\left(\Psi^{(1)}(\eta_{t}\phi)+\Psi^{(1)}((1-\eta_{t})\phi)\right)\phi^{2}\mathbf{z}_{t}\mathbf{z}_{t}'+\nabla^{(2)}\log f(\mathbf{u})
\end{eqnarray*}
In Appendix B we provide the analytical expressions for the gradient, $\nabla^{(1)}\log f(\mathbf{u})$, and the Hessian, $\nabla^{(2)}\log f(\mathbf{u})$, of the log-prior, for the different choices of the prior discussed in Section \ref{Sec_Bayesian}.
After an initial learning period the posterior mode in the various spaces of different dimension is reached with a certain tolerance value and the log-acceptance rate of the jump move is close to zero as effects of the choice of the proposal parameters.

Finally it should be noticed that the approximation $\tilde{\mathbf{u}}^{(j)}$ of the posterior mode in the jump step is also used in the M.-H. steps within the Gibbs sampler for simulating the parameters, thus no further computational burden is required for this approximation.

\section{Simulation Results}\label{Sec_SimRes}
In this section we study, through some simulation experiments, the efficiency
of the proposals for both within and between models moves. In order to check if the chain mix well we estimate the acceptance rate of the M.H. steps, the effective sample size and the convergence diagnostic statistics.

In a first set of experiments we verify the efficiency, in terms of mixing, of the proposals of the Metropolis-Hastings step of the moves for a given dimension and for the different choices of the prior distribution and parameter settings.

In each experiment we proceed as follows. For each model order $k$ and parameter setting
$\boldsymbol{\theta}'=(\boldsymbol{\alpha}',\phi)$, we generate 50 independent random samples with size $n=300$ from a $BAR(k)$. On each dataset we iterate the MCMC algorithm, defined in the previous sections, for $N=10000$ times. We consider dataset of $n=300$ observations because this is the typical sample size of the applications to real data we have in mind and that will be presented in the next section. This allows us to estimate the magnitude of the parameter estimation errors (root mean square error) associated to the inference procedure for the chosen sample size.

As regard to the values of the parameters we consider two scenarios. In the first one we
set the precision parameter $\phi=20$, which corresponds, as highlighted in Fig. \ref{Fig_BARk}, to the case of high variability in the data. We will call low precision data all the dataset simulated with this value of the precision parameter. We expect that the parameter estimation is more challenging in this context.

In the second scenario, which is more frequent in the applications we will consider in the next section, we set the precision parameter $\phi=100$ that is a higher value than the ones considered in the first scenario. In this case the data exhibits less variability (see Fig. \ref{Fig_BARk}). We will call high precision data all the dataset simulated with this precision value.

For each scenario we consider different values for the parameters and different autoregressive orders. The resulting parameter settings, which have been used in the MCMC experiments, are given in Tab. \ref{Tab_Settings}.
We also evaluate the efficiency of the MCMC algorithm for different choice of the prior distribution and consider the truncated normal distribution with parameters $\boldsymbol{\nu}=(k+2)^{-1}\boldsymbol{\iota}$ and $\Upsilon=100I$ and the
modified truncated normal distribution with parameters $\boldsymbol{\nu}=(k+2)^{-1}\boldsymbol{\iota}$, $\Upsilon=100I$ and $\kappa=10$ and the Beta-type distribution with parameters $\boldsymbol{\nu}=(k+1,\ldots,k+1)$ and $\boldsymbol{\gamma}=(k+2,\ldots,k+2)$.

\begin{table}[t]
\centering
\begin{tabular}{|c|c|c|}
\hline
$k$&$\boldsymbol{\theta}'=(\boldsymbol{\alpha}',\phi)$ &$\boldsymbol{\theta}'=(\boldsymbol{\alpha}',\phi)$\\
\hline
\hline
1&$(0.32,0.5, 20)$                                  &$(0.32,0.5, 100)$                \\
\hline
2&$(0.32,0.5, 0.1, 20)$                             &$(0.32,0.5, 0.1, 100)$           \\
\hline
3&$(0.32,0.5, 0.1, 0.03, 20)$                       &$(0.32,0.5, 0.1, 0.03, 100)$     \\
\hline
4&$(0.32,0.4, 0.1, 0.03,0.1, 20)$                   &$(0.32,0.4, 0.1, 0.03, 0.1, 100)$\\
\hline
\end{tabular}
\caption{The parameter settings employed in the MCMC experiments. \textit{First column}: the autoregressive order. \textit{Second and third columns}: conditional mean parameters for low precision data (i.e. $\phi=20$) and high precision data (i.e. $\phi=100$).}\label{Tab_Settings}
\end{table}

\begin{figure}[t]
\centering
\includegraphics[width=200pt]{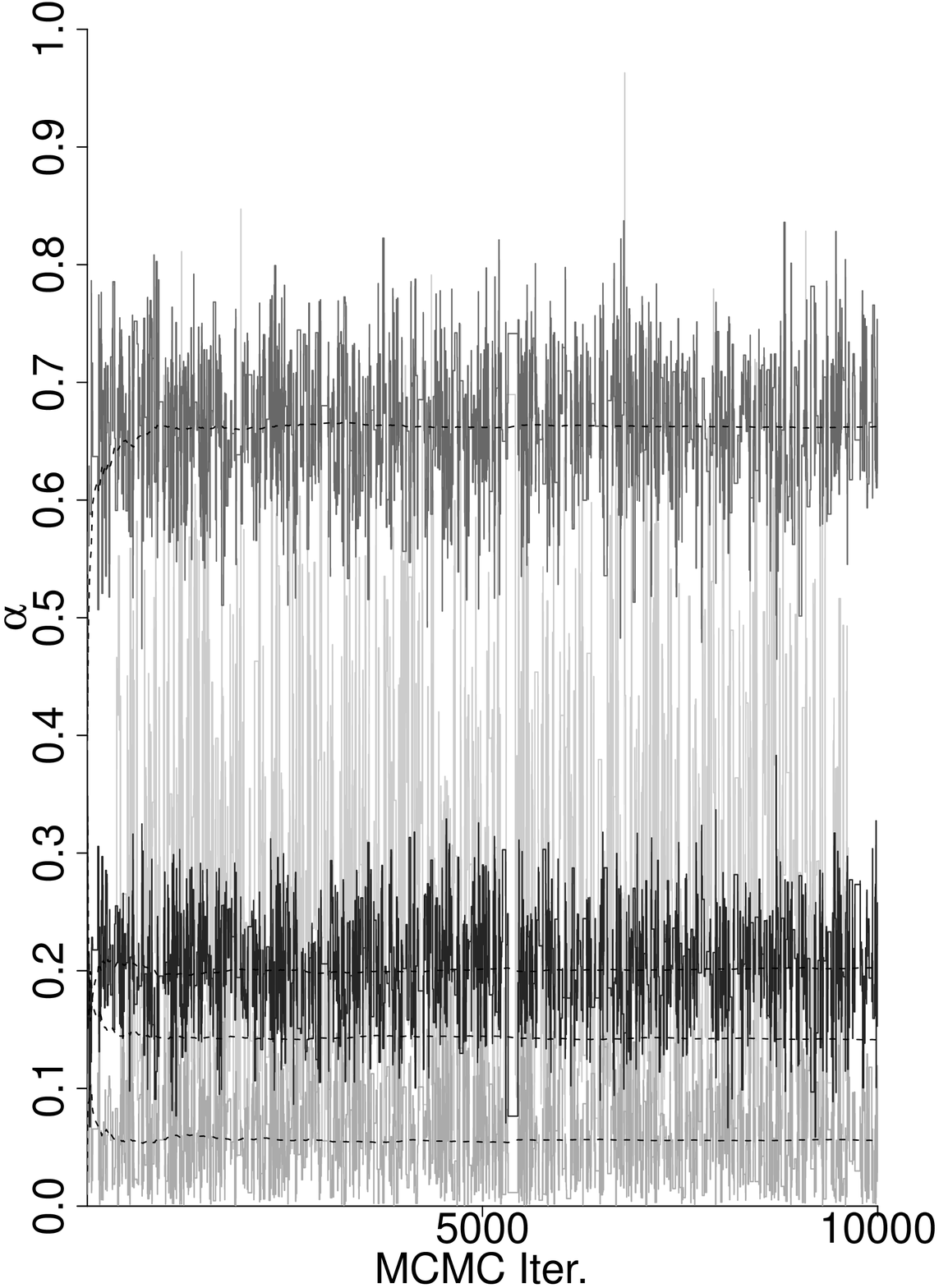}
\includegraphics[width=200pt]{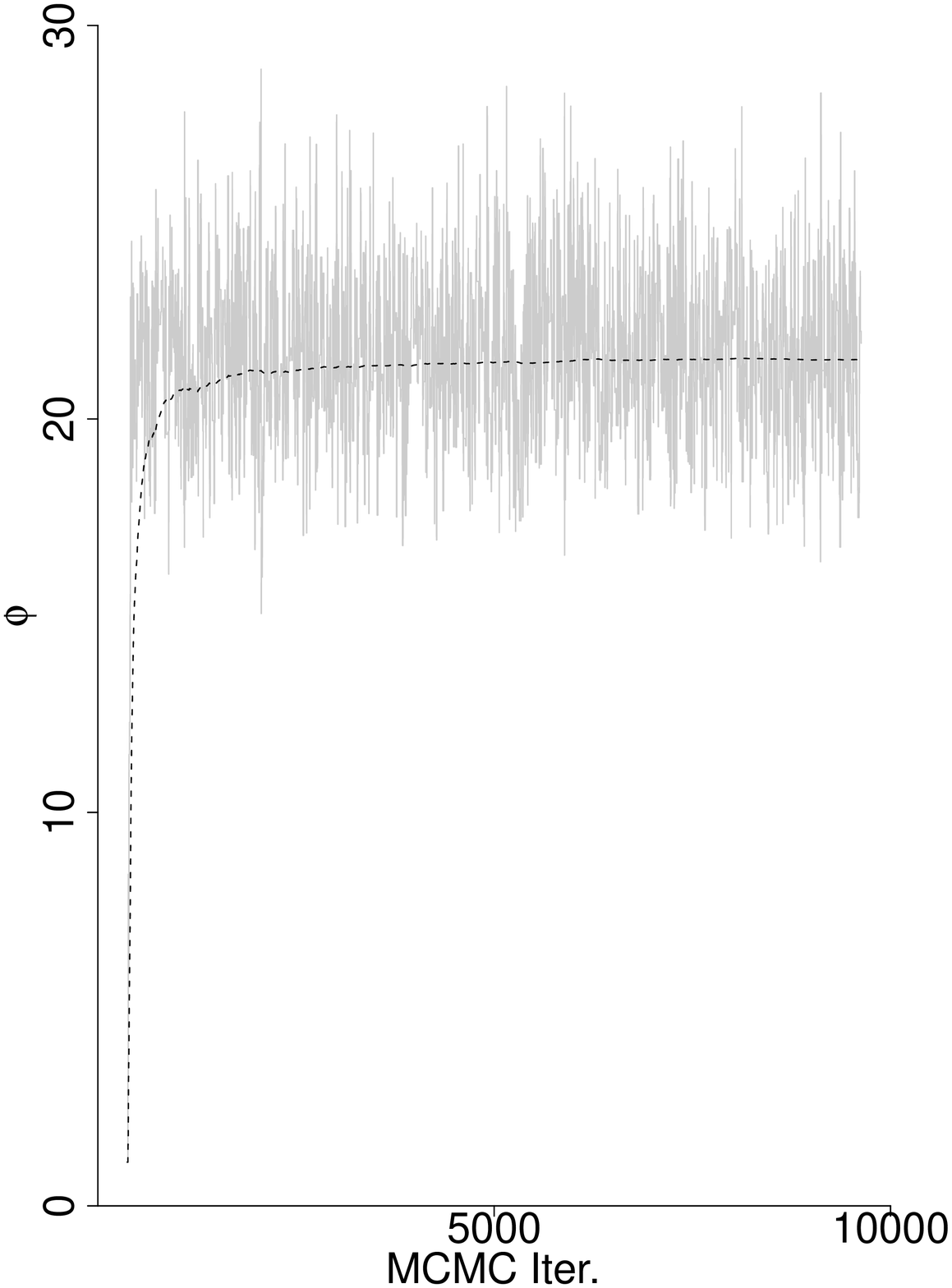}
\caption{MCMC raw output (\textit{gray solid lines}) and progressive averages (\textit{dashed lines}) for $\boldsymbol{\alpha}$ (\textit{left chart}) and $\phi$ (\textit{right chart}) estimated on a dataset of $n=300$ observations simulated from a BAR(3) process with parameters $\boldsymbol{\alpha}=(0.32,0.5, 0.1, 0.03)'$ and $\phi=20$. The initial value of the Gibbs sampler and the burn-in sample are included in the MCMC sample in order to show the convergence of the MCMC progressive averages to the true values of the parameters.} \label{Fig_BARMCMC}
\end{figure}

Fig. \ref{Fig_BARMCMC} shows typical raw output and progressive averages of $N=10000$ iterations of the MCMC chain for $\boldsymbol{\alpha}$ (left chart) and $\phi$ (right chart). In these figures, the initial value of the Gibbs sampler and the burn-in sample are included in order to show the convergence of the MCMC progressive averages to the true values of the parameters. In this example, we considered a sample of $n=300$ observations simulated from a BAR(3) model with parameters $\boldsymbol{\alpha}=(0.32,0.5, 0.1, 0.03)'$ and $\phi=20$.

For each set of 50 independent MCMC experiments of length $N$ we estimate the RMSE, the average acceptance probability (ACC) of the M.-H. within Gibbs steps and the following quantity
\begin{equation}
ESS=N\left(1+\sum_{t=1}^{\infty}\hbox{corr}(\theta^{(0)},\theta^{(t)})\right)^{-1}
\end{equation}
that is the effective sample size (ESS) of the MCMC sample for a parameter $\theta$.
In order to evaluate the RMSE and the other statics we consider $N=10000$ MCMC iterations and discard the output of the first 1000 iterations.

In the various experiments we will also evaluate the number of MCMC iterations which are necessary to our MCMC algorithm to reach convergence. It should be noticed that the choice of the number of MCMC iterations is still an open issue (see \citet{RobCas04} and \citet{GuiMenRob98}). In this work we combine a graphical inspection of the progressive averages of the parameter posterior distribution with the evaluation  of a convergence criteria based on the Kolmogorov-Smirnov (KS) test (see \citet{RobCas04}, Ch. 12). See also \cite{BroGiuPhi03} for an alternative use of the KS statistics for detecting convergence in RJMCMC chains. In the simulation experiments, for each component $\theta$ of the parameter vector $\boldsymbol{\theta}$, we split the associated MCMC sample $\theta^{(j)}$, $j=1,\ldots,N$
in two subsamples $\theta_{1}^{(j)}$ and $\theta_{2}^{(j)}$ with $j=1,\ldots,M$
and evaluate
\begin{equation}
KS=\frac{1}{M}\underset{\eta}{\sup}\left|\sum_{j=1}^{M}\mathbb{I}_{(0,\eta)}(\theta_{1}^{(jG)})-\sum_{j=1}^{M}\mathbb{I}_{(0,\eta)}(\theta_{2}^{(jG)})\right|
\end{equation}
where $G$ is the batch size. The use of batches is necessary in order to obtain quasi-independent samples. The independence of the samples is one of the assumption to have a known limiting distribution for the KS statistics. For each experiment we will show the average p-value of the KS statistics over the vector of parameters and over the last 100 iterations of the MCMC chain.

We summarize the results of the MCMC experiments in Tab. \ref{Tab_MCMC}. The results for the Beta-type prior and for the modified truncated-Gaussian are similar, thus the choose to present the results for modified truncated-Gaussian prior. From Tab. \ref{Tab_MCMC} we note that the RMSE of the autoregressive coefficients and of the $\phi$ parameter tends to grow with the model order. This suggests that the higher the number of parameters to be estimated the lower the precision of the estimates. Besides, the precision of the $\phi$ estimate is influenced by the number of parameters, as it decreases as the model order increases.

\begin{table}[p]
\centering
\begin{tabular}{|c|cccccc|c|c|c|}
\hline
\multicolumn{10}{|c|}{Low precision data}\\
\hline
\hline
k&\multicolumn{6}{|c|}{Estimated RMSE}                        & $ACC$ & $ESS$ & $KS$\\
 &$\alpha_{0}$&$\alpha_{1}$&$\alpha_{2}$&$\alpha_{3}$&$\alpha_{4}$&$\phi$&               &       &     \\
\hline
\hline
&\multicolumn{9}{|c|}{\textit{Truncated-Gaussian Prior} ($\boldsymbol{\nu}=(k+2)^{-1}\boldsymbol{\iota}$, $\Upsilon=100I$)}\\
\hline
1&0.032&0.058&     &      &     &0.376 &  0.176    & 704   &0.534 \\
\hline
2&0.033&0.043&0.023&     &      &0.996 &  0.172    & 630   &0.552 \\
\hline
3&0.087&0.094&0.026&0.051&      &2.092 &  0.163    & 584   &0.523 \\
\hline
4&0.041&0.011&0.019&0.075&0.032 &3.727 &  0.155    & 538   &0.541 \\
\hline
\hline
&\multicolumn{9}{|c|}{\textit{Modified Truncated-Gaussian Prior} ($\boldsymbol{\nu}=(k+2)^{-1}\boldsymbol{\iota}$, $\Upsilon=100I$, $\kappa=10$)}\\
\hline
1&0.033&0.051&     &     &     &0.392 &   0.181    & 1013   &0.556 \\
\hline
2&0.032&0.055&0.021&     &     &0.916 &   0.183    & 853   &0.563 \\
\hline
3&0.015&0.077&0.018&0.059&     &1.701 &   0.192    & 783   &0.574 \\
\hline
4&0.030&0.023&0.013&0.059&0.034&2.564 &   0.189    & 740   &0.539 \\
\hline
\multicolumn{10}{c}{\vspace{1pt}}\\
\hline
\multicolumn{10}{|c|}{High precision data}\\
\hline
\hline
k&\multicolumn{6}{|c|}{Estimated RMSE}                        & $ACC$ & $ESS$ & $KS$\\
 &$\alpha_{0}$&$\alpha_{1}$&$\alpha_{2}$&$\alpha_{3}$&$\alpha_{4}$&$\phi$&               &       &     \\
\hline
\hline
&\multicolumn{9}{|c|}{\textit{Truncated-Gaussian Prior} ($\boldsymbol{\nu}=(k+2)^{-1}\boldsymbol{\iota}$, $\Upsilon=100I$)}\\
\hline
1&0.011&0.018&     &      &     & 0.964&  0.392    &  923  &0.513\\
\hline
2&0.029&0.047&0.031&      &     & 1.815&  0.402    &  827  &0.546\\
\hline
3&0.038&0.071&0.032&0.002&      & 3.122&  0.422    &  798  &0.593\\
\hline
4&0.021&0.038&0.037&0.007&0.029 & 6.430&  0.538    &  778  &0.511\\
\hline
\hline
&\multicolumn{9}{|c|}{\textit{Modified Truncated-Gaussian Prior} ($\boldsymbol{\nu}=(k+2)^{-1}\boldsymbol{\iota}$, $\Upsilon=100I$, $\kappa=10$)}\\
\hline
1&0.017&0.021&     &     &      &0.392  &  0.403     & 1198  &0.511\\
\hline
2&0.020&0.028&0.001&     &      &1.101  &  0.420     & 1012   &0.534\\
\hline
3&0.031&0.063&0.003&0.002&      &1.539  &  0.428     & 941   &0.529\\
\hline
4&0.029&0.033&0.033&0.001&0.018 &3.955  &  0.509     & 830   &0.542\\
\hline
\end{tabular}
\caption{Estimation results for different model orders, parameter settings and prior distributions. The results are averages over a set of 50 independent MCMC experiments on 50 independent dataset of $n=300$ observations. On each dataset we run the proposed MCMC algorithm for $N=10000$ iterations and then discard the first 1000 iterations in order to estimate the parameters.
In each panel: model order (\textit{first column}); estimated root mean square error (RMSE) (\textit{second column}) for each parameter; average acceptance rate (ACC) and the effective sample size (ESS) (\textit{third and fourth columns}) averaged over the two M.-H. chain in the Gibbs sampler; KS convergence diagnostic statistics (\textit{last column}) with batch size $G=50$ (average over the last 100 iterations).}\label{Tab_MCMC}
\end{table}

Moreover, the RMSEs of the high precision data results are generally lower than RMSEs of low precision data, with the exception of the $\phi$ parameter that seems to behave differently from the autoregressive coefficients. However, in order to correctly evaluate the outcomes, we need to compare the RMSE to the value of the $\phi$ parameter. For example, for the Truncated-Gaussian prior the estimated RMSE of $\phi$ for the low precision $BAR(4)$ is equal to $3.727$ and the corresponding high precision RMSE is equal to $6.430$. However, dividing this two values by the related $\phi$ value ($20$ and $100$ respectively) we get percentages of  $18.635$ and $6.43$, showing a better performance for high precision results, as expected.

The third column of Tab. \ref{Tab_MCMC} lists the average acceptance rate of the parameters. As it is clear from the table, data with higher precision show a noticeably improvement in terms of ACC compared to data with lower precision. However, we need to point out that Tab. 2 displays the average acceptance rate and the improvement is due mainly to the higher ACC value of the $\phi$ parameter, since the ACC of the autoregressive coefficients is very similar for both type of data.

The ESS is the number of effectively independent draws from the posterior distribution and it is a measure of the mixing of the Markov chain. In our case, this measure decreases and denotes a worse mixing as the order of the process increases. Moreover, as we see in the fourth column of Tab. 2, ESS values are higher in high precision data, showing a better mixing than in low precision data.

The average p-values of the KS statistic, in the last column of Tab. \ref{Tab_MCMC}, take values close to 0.5 for the different precision of the data and the different model orders, suggesting the acceptance of the null hypothesis that the subsamples associated to the Markov chain have the same distribution, guaranteeing convergence. Besides, since Tab. 2 displays only the average p-values, we underline that over the last 100 iterations of the Markov chains the KS p-values improved, getting closer to 1.

Comparing now the outcomes related to the priors used, the Modified Truncated-Gaussian prior performs better than the Truncated-Gaussian prior for both the low and high precision data. As we notice from Tab. \ref{Tab_MCMC}, the RMSEs of the autoregressive coefficients are quite similar, but the RMSE of the $\phi$ parameter calculated with the modified prior is sensibly lower. Moreover, the ACC and the ESS values of the modified prior are slightly higher, denoting a better mixing of the Markov chains. The p-values of the KS statistic indicate convergence to the target distribution, being close to the value of 0.5.

\begin{figure}[p]
\centering
\includegraphics[width=320pt]{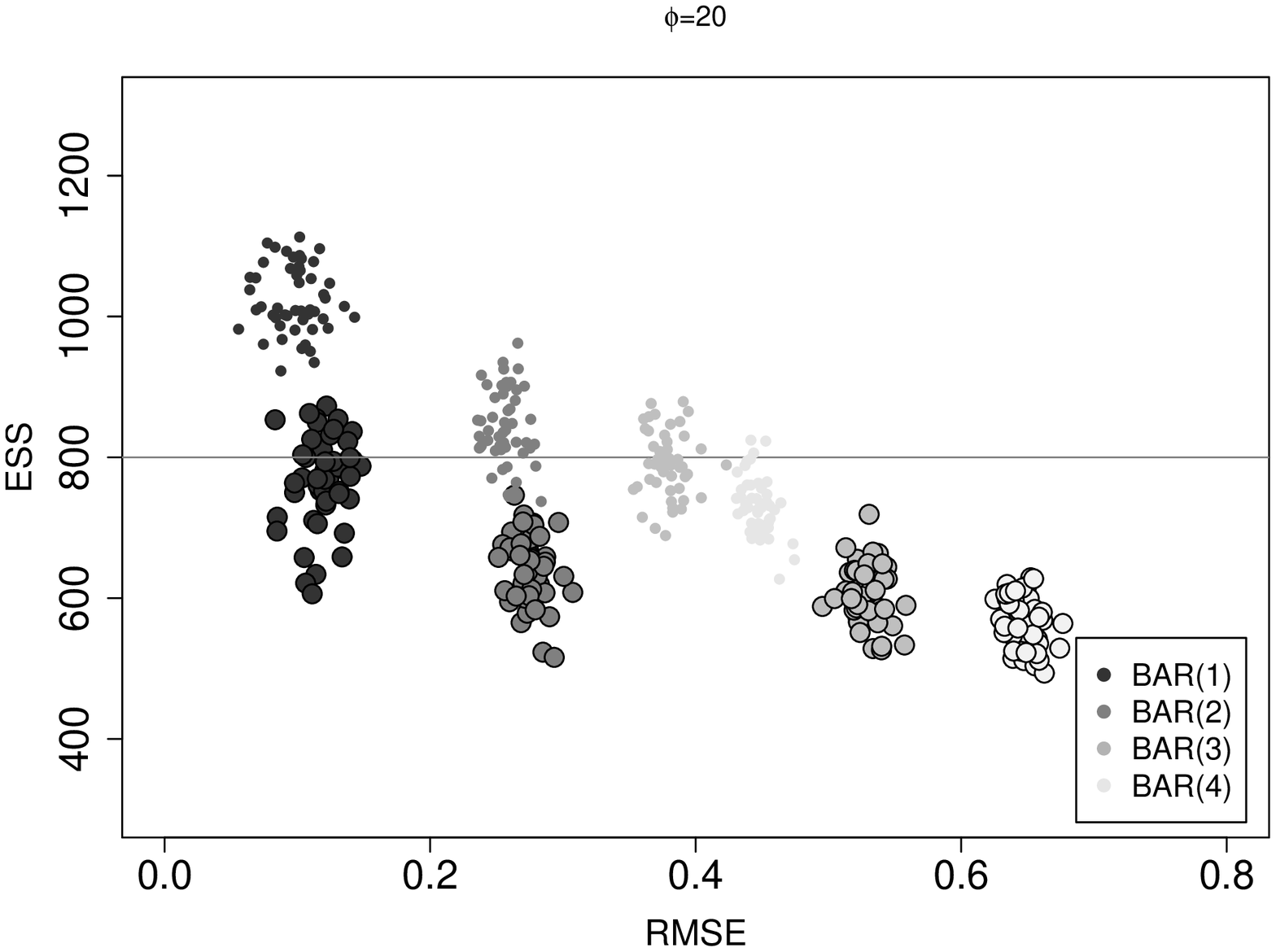}
\includegraphics[width=320pt]{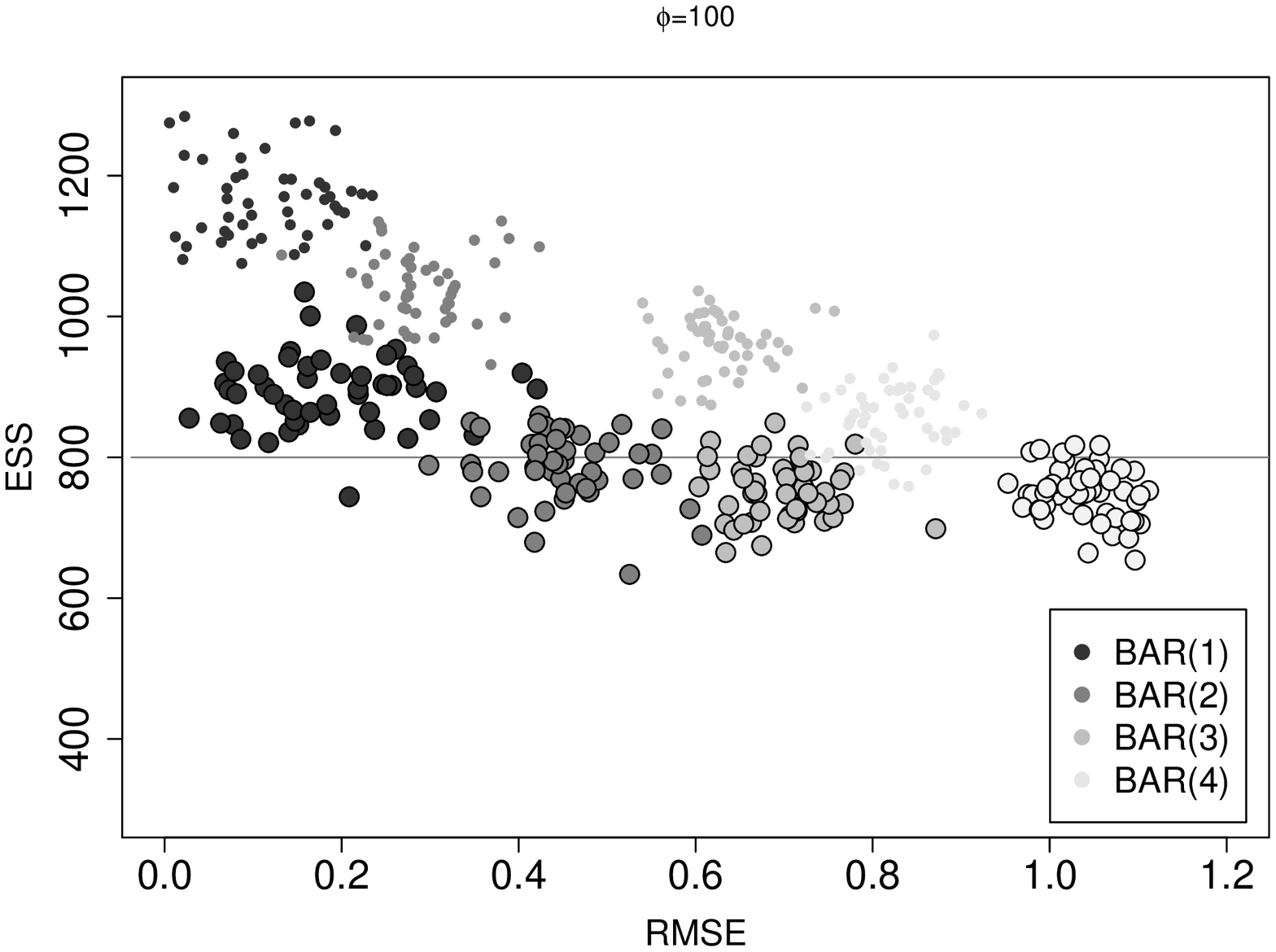}
\caption{ESS and RMSE of the 50 MCMC simulation experiments (averages over the parameter vector) for different autoregressive order of the $BAR(k)$ (different gray levels of the circles) and different choice of the prior (bigger circles for the truncated Gaussian and smaller circles for the modified Gaussian prior). In the rows the ESS and the RMSE statistics for two different values of $\phi$.} \label{Fig_BARMCMC1}
\end{figure}

Fig. \ref{Fig_BARMCMC1} gives a more detailed description of the behavior of the RMSE and of the ESS in the MCMC experiments.  Fig. \ref{Fig_BARMCMC1} illustrates the parameters RMSE on the $x$-axis and the ESS of the two M.H. steps on the $y$-axis, for the 50 MCMC simulation experiments. Large circles denote the truncated-Gaussian prior results, while small circles denote the modified truncated-Gaussian prior results. The different colors of the circles indicate the different order of the $BAR(k)$. The top plot displays low precision results ($\phi=20$) and the bottom plot displays high precision results ($\phi=100$).

For both high and low precision data, the RMSE grows with the order of the model, while the ESS decreases. Moreover, the modified prior gives more efficient estimates, since the small circles in Fig. \ref{Fig_BARMCMC1} show an improvement as in the RMSE as in the ESS.

The RMSE calculated with low precision data ($\phi=20$) is lower that the RMSE calculated with high precision data ($\phi=100$), but it is due to the fact that the circles depicted in Fig. \ref{Fig_BARMCMC1} are averages over all the parameters, including the parameter $\phi$. Therefore, as illustrated above, in order to compare the outputs it is necessary to consider the different value of the parameters.

\begin{figure}[p]
\centering
\includegraphics[width=320pt]{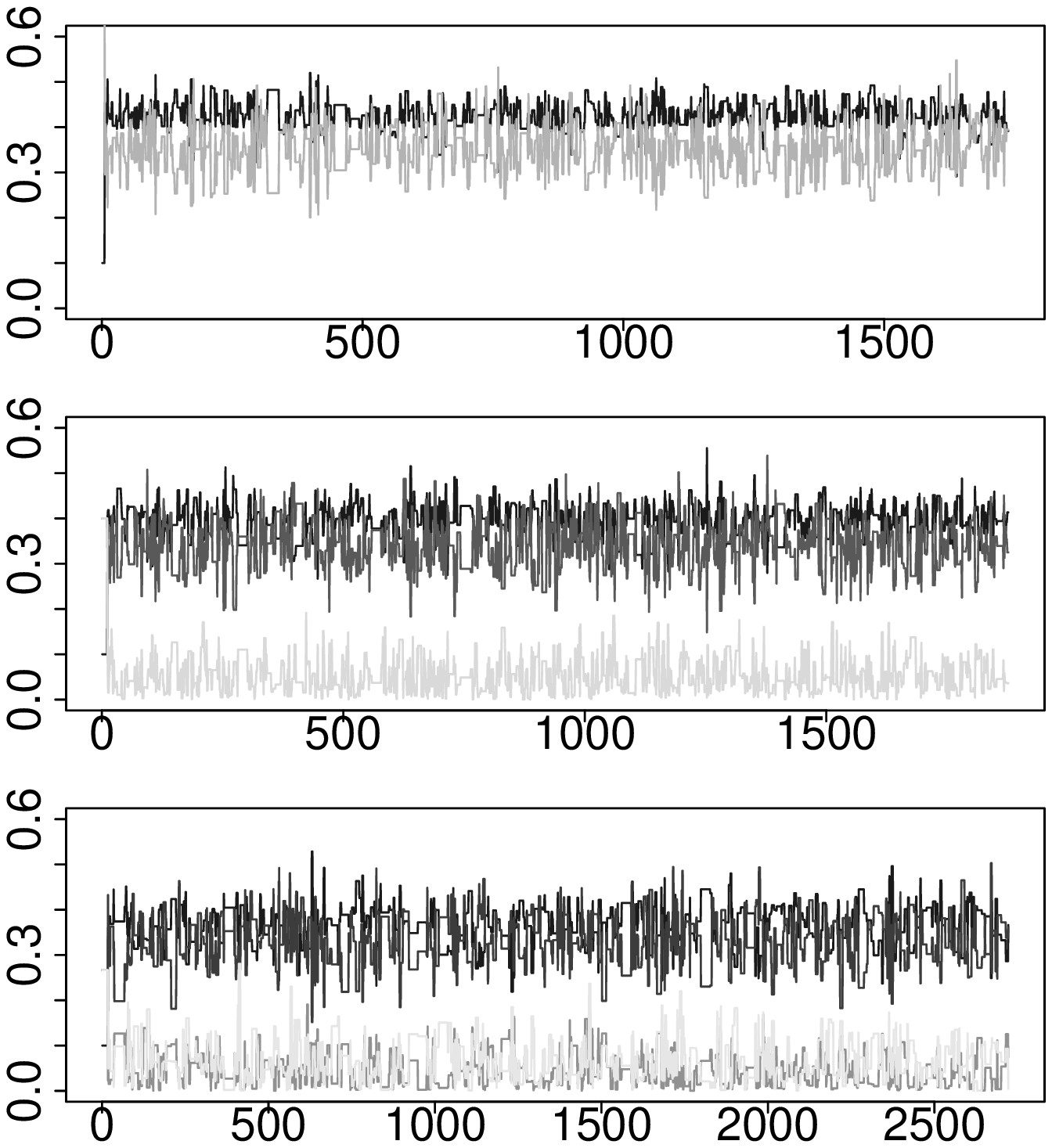}
\includegraphics[width=320pt]{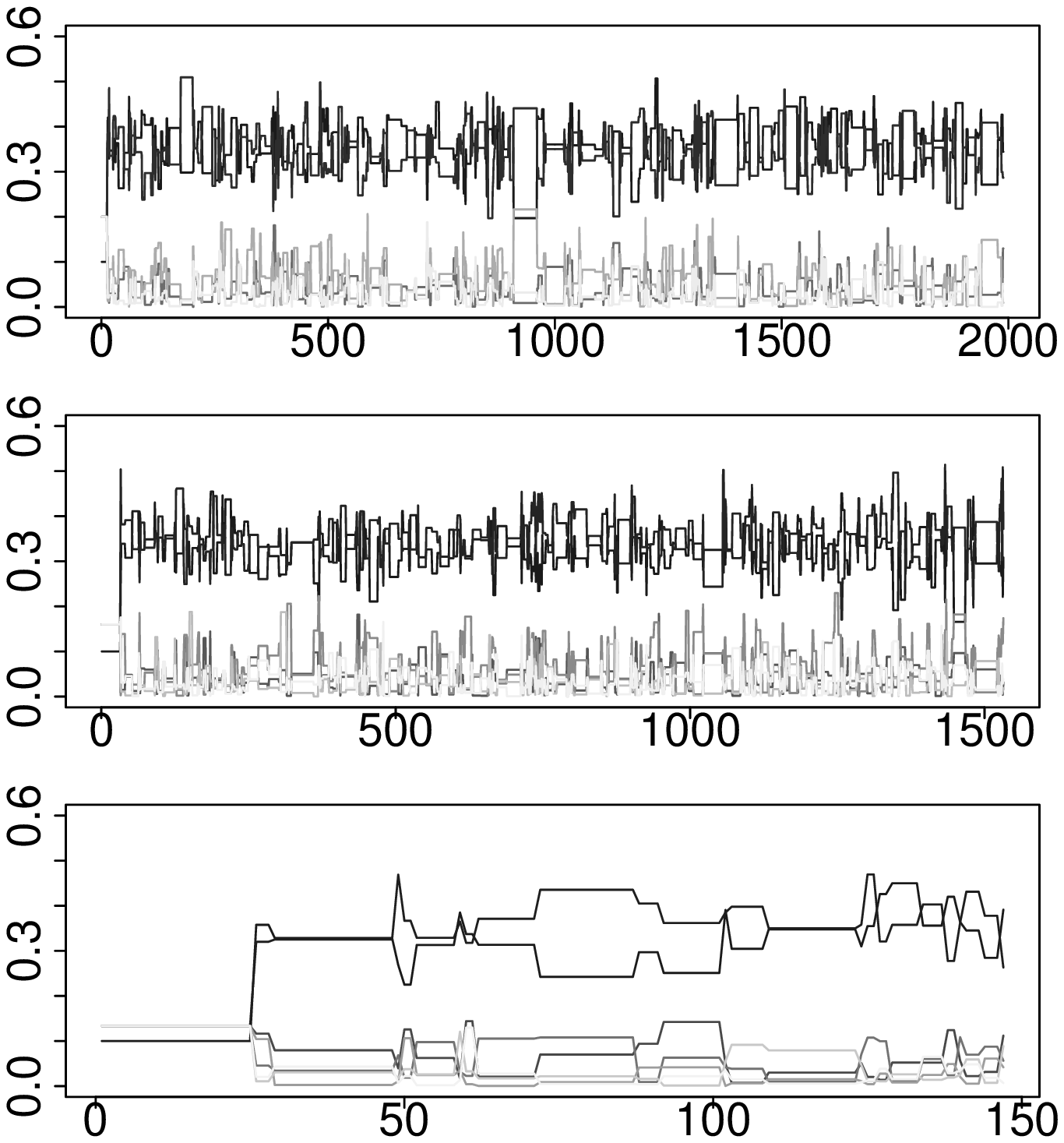}
\caption{RJMCMC model selection when the true model is a $BAR(3)$ with $\boldsymbol{\alpha}=(0.37,0.4,0.1,0.03)'$
and $\phi=100$. In each graph, the raw output of the RJMCMC chain (\textit{gray lines}) for $\boldsymbol{\alpha}\in\Delta_{k+1}$ in spaces with different dimensions $(k+1)$ with $k=1,\ldots,6$.} \label{Fig_BARRJMCMC1}
\end{figure}

\begin{figure}[h]
\centering
\includegraphics[width=320pt]{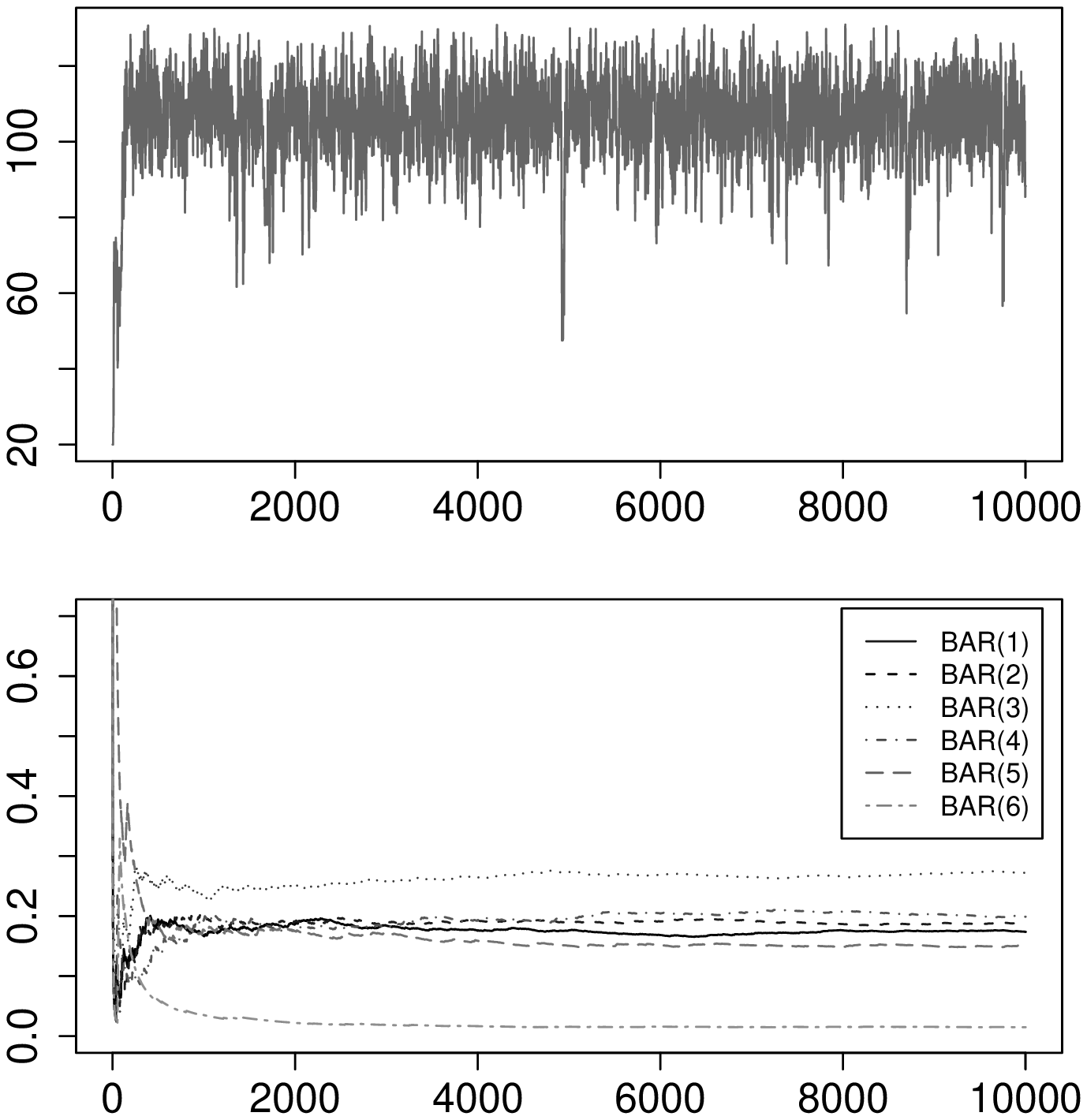}
\caption{RJMCMC raw output for $\phi$ (\textit{upper chart}),
the progressive estimates of the model probabilities (\textit{bottom chart}), when the true model is a $BAR(3)$
with $\boldsymbol{\alpha}=(0.37,0.4,0.1,0.03)'$ and $\phi=100$.} \label{Fig_BARRJMCMC2}
\end{figure}

Furthermore, we note that for the different model orders high precision results appear to be more scattered, while low precision outcomes appear to be more concentrated. Since this behavior is associated to a higher ESS value, this suggest that the Markov chains explore more freely the parameter space and that they are characterized by a better mixing.

The aim of the second set of experiments is to study the relation between the size of the sample of observations and the model posterior probability.
We consider a dataset of 500 observations simulated from a BAR(k) with $k=3$ and parameter values $(\boldsymbol{\alpha}',\phi)=(0.37,0.4, 0.1, 0.03, 100)$. Then in order to estimate both the parameters and the autoregressive order we assume a modified truncated-Gaussian prior with parameters $\boldsymbol{\nu}=(k+2)^{-1}\boldsymbol{\iota}$, $\Upsilon=100I$ and apply the RJMCMC algorithm presented in the previous section to subsamples of different size (from 100 to 500 observations).

A typical output of a RJMCMC chain for a dataset of 300 observations is given in Fig. \ref{Fig_BARRJMCMC1}. Each chart shows the RJMCMC iterations for each subspace. The iterations of the chain for the precision parameter $\phi$ are given in the upper chart of Fig. \ref{Fig_BARRJMCMC2}. In the bottom chart of the same figure we show the model posterior probabilities.
\begin{table}[t]
\centering
\begin{tabular}{|c|ccccc|}
\hline
    \multicolumn{6}{|c|}{Estimated Model Probabilities}\\
\hline
\backslashbox{k}{n} & 100& 200&300 &400&500   \\
\hline
 1  &0.06667&0.09143&0.10291&0.09806&0.00129\\
 2  &0.19208&0.38611&0.29345&0.19384&0.17828\\
 3  &0.18939&0.43252&0.35730&0.65473&0.75391\\
 4  &0.13245&0.03387&0.24522&0.05255&0.04449\\
 5  &0.10674&0.02291&0.00002&0.00011&0.00018\\
 6  &0.06682&0.01910&0.00001&0.00033&0.00013\\
 7  &0.05682&0.00691&0.00001&0.00001&0.00663\\
 8  &0.04682&0.00149&0.00001&0.00001&0.00012\\
 9  &0.04283&0.00139&0.00001&0.00001&0.00018\\
 10 &0.03891&0.00109&0.00042&0.00001&0.00001\\
 11 &0.01488&0.00080&0.00014&0.00001&0.00001\\
 12 &0.01455&0.00055&0.00022&0.00001&0.00001\\
 13 &0.01109&0.00062&0.00012&0.00001&0.00001\\
 14 &0.01078&0.00027&0.00011&0.00001&0.00001\\
 15 &0.00917&0.00092&0.00005&0.00030&0.00023\\
\hline
 Mode  &2      &3      &3      &3         &3      \\
 Mean  &4.80121&2.65058&2.75522&2.66801   &2.85226\\
 s.d.  &3.1253&1.2256  &0.9854 &0.76153   &0.70101\\
\hline
\end{tabular}
\caption{Relation between sample size $n$ (\textit{first row}) and model order posterior (\textit{columns from one to six}) for $k_{\max}=15$ when data are simulated from a BAR(3) with $(\boldsymbol{\alpha}',\phi)=(0.37,0.4, 0.1, 0.03, 100)$. We assume a modified truncated-Gaussian prior with parameters $\boldsymbol{\nu}=(k+2)^{-1}\boldsymbol{\iota}$, $\Upsilon=100I$. Approximation of the model order posterior and of its mode, mean and
standard deviation (\textit{last three rows}) is based on 100,000 RJMCMC iterations.}\label{Tab_RJMCMC}
\end{table}

The results of the model selection procedure for the different sample size are given in Tab. \ref{Tab_RJMCMC}. All the estimates are based on 100,000 iterations of the RJMCMC chain. In the specific example at hand, the last three rows in the table show two interesting results. First the estimated model order is not correct ($\hat{k}=2$ with $k=3$) for a sample of 100 observations. Secondly the standard deviation is $3.1253$ for $n=100$ and decreases for increasing $n$, which means that the model posterior distribution is less concentrated around the mode for the samples of smaller size.

It should be noticed that, while small-sample bias in the model order estimates has also been observed in RJMCMC-based model selection procedures for other non-Gaussian autoregressive processes, such as the integer valued ARMA processes (see \cite{EncNeaSub09}), what we found new for the Beta autoregressive, is the high dispersion of the model order estimates in small samples.

Both of these results may contribute to explain the shape of the model posterior distribution that we will obtain in the applications to real dataset presented in Section \ref{Sec_Emp}.

\section{Empirical Application}\label{Sec_Emp}
\begin{figure}[t]
\centering
\includegraphics[width=320pt]{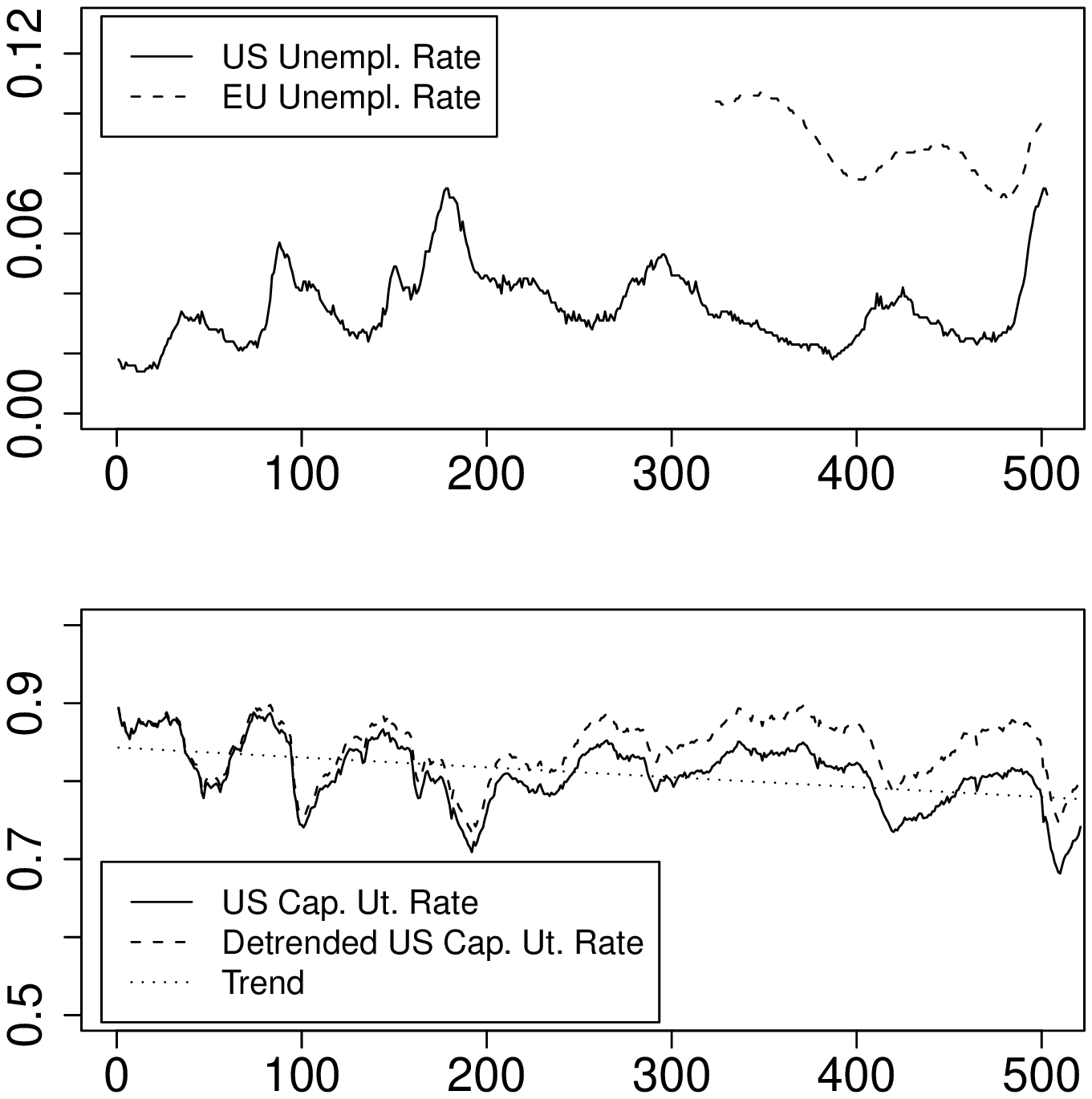}
\caption{\textit{Up}: US (\textit{solid line}) and EU (\textit{dashed line}) unemployment rates at a monthly frequency. \textit{Bottom}: US capacity utilization rate (\textit{solid line}), its estimated linear trend (\textit{dotted line}) with intercept $\hat{\gamma_{0}}=0.843$ and slope $\hat{\gamma}_{1}=-0.066$ and the detrended capacity rate (\textit{dashed line}).} \label{Fig_Unempl}
\end{figure}
\subsection{Unemployment Data}
The dynamics of the unemployment rate is certainly one of the most studied in the economics time series analysis (see \cite{Bea94} and \cite{Nic97} for an early review on unemployment models). Modeling and forecasting the unemployment rate still represent some of the most challenging issues in econometrics. See for example \cite{Nef84}, \cite{MonZarTsaTia98} and \cite{KooPot99} for some recent advances on nonlinear models. The unemployment rate is usually characterized by relatively brief periods of rapid economic contraction and by relatively long periods of slow expansion. In this work we do not model the asymmetric behavior of time series but focus instead on another fundamental feature of this variable, that is the unemployment rate is naturally defined on a bounded interval, let us say the $(0,1)$ interval. Data transformation is usually applied to have data on the real interval (see \cite{Wal87}). Recently \citet{RocCri09} suggested instead Beta autoregressive models for modeling the unemployment rate.

Here we consider two interesting dataset (see Fig. \ref{Fig_Unempl}). The first one is the US unemployment rate (source: Datastream) sampled at a monthly frequency from February 1971 to December 2009. We are mainly interested in modeling the economic cycle and thus consider deseasonalized data. This dataset is quite large (467 observations) when compared to other macroeconomics dataset and is one of the most studied set of data in econometrics (see for example \cite{Nic97}). The other dataset is the deseasonalized unemployment rate of the Euro Area sampled at the monthly frequency (source: Datastream), from January 1995 to December 2009. We will consider the aggregated unemployment rate for the 13 countries area. This is another well studied dataset (see \cite{Bea94} and inference on this dataset could be challenging due to the limited amount of observations (180 observations). Moreover modeling and forecasting of this variable represents one of the most important issues for the European Central Bank and for the European institute of official statistics (Eurostat).

We assume modified Gaussian prior for the autoregressive parameters and a uniform prior for the autoregressive order. For the RJMCMC algorithm we consider $N=100,000$ iterations and a discard the first 10,000 samples in order to obtain an estimate of the parameters and of the model posterior. The estimation results are given in Tab. \ref{Tab_Unempl}. Figures \ref{Fig_BARRJMCMCUnemploy} show the model posterior probabilities for the two time series.

\begin{table}[t]
\centering
\begin{tabular}{|c|cccccccc|c|}
\hline
\multicolumn{10}{|c|}{US Unemployment Rate ($\hat{k}=2$)}\\
\hline
$\boldsymbol{\theta}_{k}$      & $\alpha_{0}$ & $\alpha_{1}$ & $\alpha_{2}$ &              &              &              &               &$\phi$& ACC\\
\hline
$\hat{\boldsymbol{\theta}_{k}}$& 0.011        & 0.547        &  0.272       &              &              &              &              & 130   & 0.392\\
\hline
$Q_{0.025}$                    & 0.006        & 0.118        &  0.007       &              &              &              &              & 126   &    \\
$Q_{0.975}$                    & 0.017        & 0.820        &  0.708       &              &              &              &              & 131   &    \\
\hline
\hline
\multicolumn{10}{|c|}{EU Unemployment Rate  ($\hat{k}=5$)}\\
\hline
$\boldsymbol{\theta}_{k}$      & $\alpha_{0}$ & $\alpha_{1}$ & $\alpha_{2}$ & $\alpha_{3}$ & $\alpha_{4}$ & $\alpha_{5}$  &              & $\phi$&ACC \\
\hline
$\hat{\boldsymbol{\theta}_{k}}$& 0.029       &  0.544       &  0.076       & 0.010        &  0.024       &  0.011       &              & 128   & 0.278 \\
\hline
$Q_{0.025}$                    & 0.024       &  0.140       &  0.001       & 0.004        &  0.010       &  0.004       &              & 123   &      \\
$Q_{0.975}$                    & 0.034       &  0.572       &  0.096       & 0.016        &  0.029       &  0.016       &              & 132   &      \\
\hline
\hline
\multicolumn{10}{|c|}{US Capacity Utilization Rate  ($\hat{k}=6$)}\\
\hline
$\boldsymbol{\theta}_{k}$      & $\alpha_{0}$ & $\alpha_{1}$ & $\alpha_{2}$ & $\alpha_{3}$ & $\alpha_{4}$ & $\alpha_{5}$  & $\alpha_{6}$ & $\phi$&ACC \\
\hline
$\hat{\boldsymbol{\theta}_{k}}$&  0.397      &  0.196      & 0.075       &   0.065     &  0.045      &   0.053     &  0.049       & 130   & 0.398 \\
\hline
$Q_{0.025}$                    &  0.390      &  0.189      & 0.058       &   0.058     &  0.040      &   0.049     &  0.030      & 126   &      \\
$Q_{0.975}$                    &  0.408      &  0.207      & 0.079       &   0.074     &  0.067      &   0.060     &  0.051      & 130   &      \\
\hline
\end{tabular}
\caption{Estimated parameters $\boldsymbol{\theta}$ and the model order $k$ for the US (\textit{first panel}) and EU13 (\textit{second panel}) unemployment rates
and for the US capacity utilization rate (\textit{last panel}). In the last column the acceptance rate (ACC) of the trans-dimensional jump move in the RJMCMC chain.}\label{Tab_Unempl}
\end{table}

\subsection{Capacity Utilization}
We consider the capacity utilization rate. While unemployment rate deals with utilization of the labor as a production factor, the capacity utilization deals with all production factors (i.e. labor force and stock of capital) of an economic system or sector.  A detailed definition of capacity utilization and a discussion on the different ways to get a statistical measure of this quantity can be found in \cite{KleSu79}.

The capacity utilization is a relevant quantity in both economic theory (see for example \cite{BurEic96} and \cite{CooHanPre95}) and in time series macroeconometrics (see \cite{KleSu79}). In time series analysis both modelling and forecasting of this indicator are challenging issues which have a crucial role in the practice of economics analysis. In fact a decreasing capacity utilization is usually interpreted as slowdown of the aggregated demand and consequently a reduction of the inflation level. An increase of the capacity utilization reveals an expansion of the level of economic activity.

In this work we consider the capacity utilization rate series for the US sampled at a monthly frequency from January 1967 to May 2010. The series refers to all the industry sectors and is seasonally adjusted (source: Datastream).

\begin{figure}[p]
\centering
\includegraphics[width=220pt]{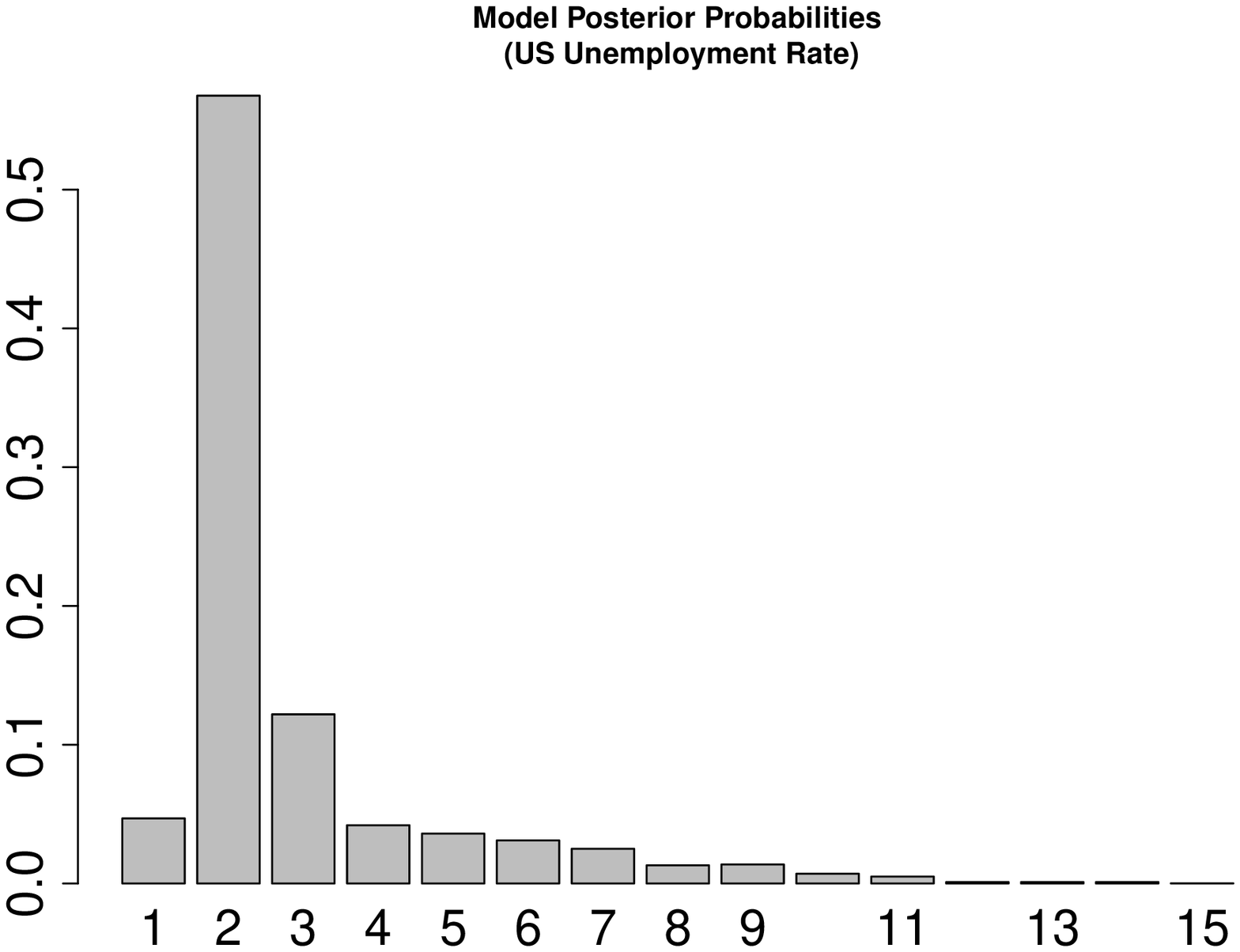}
\includegraphics[width=220pt]{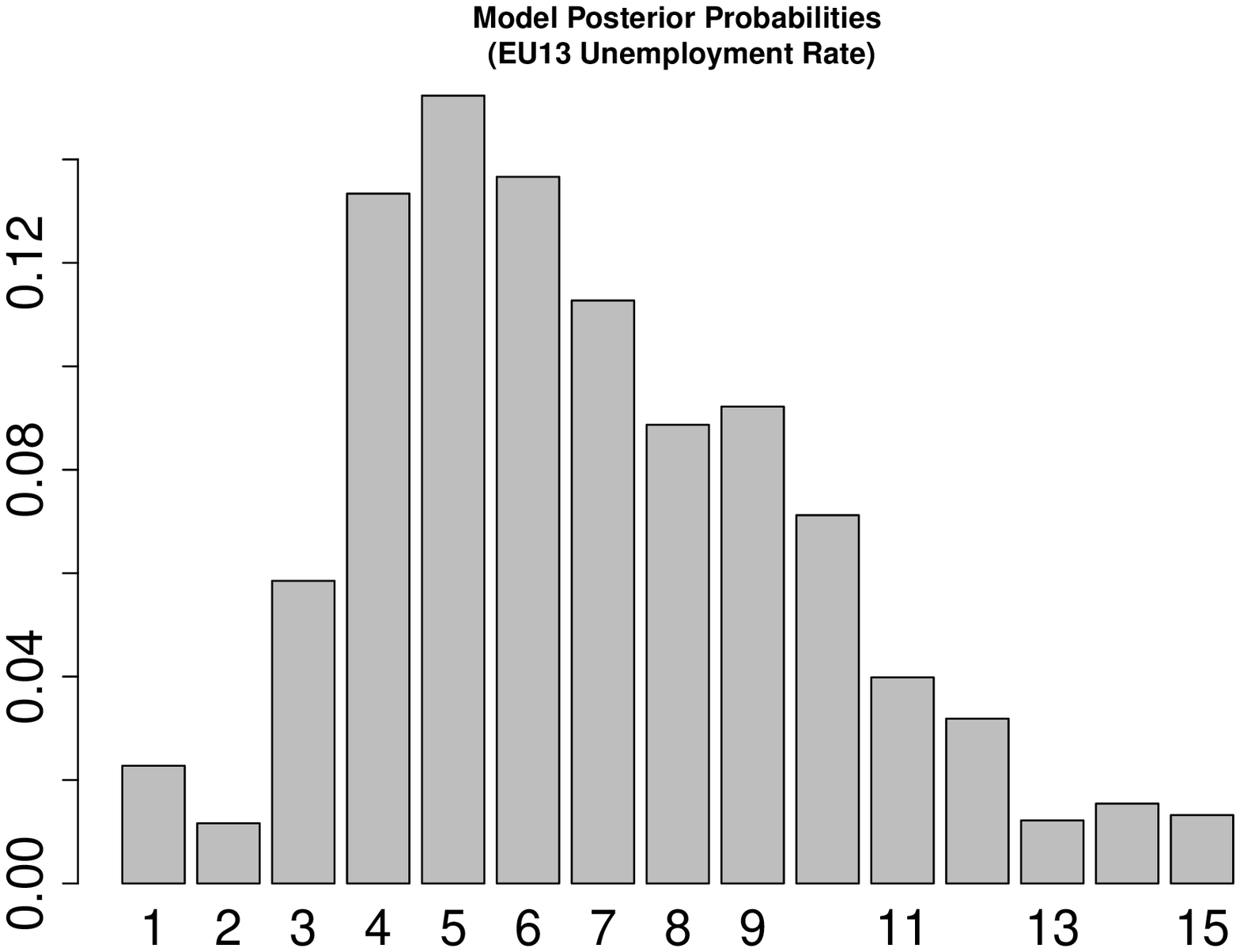}
\includegraphics[width=220pt]{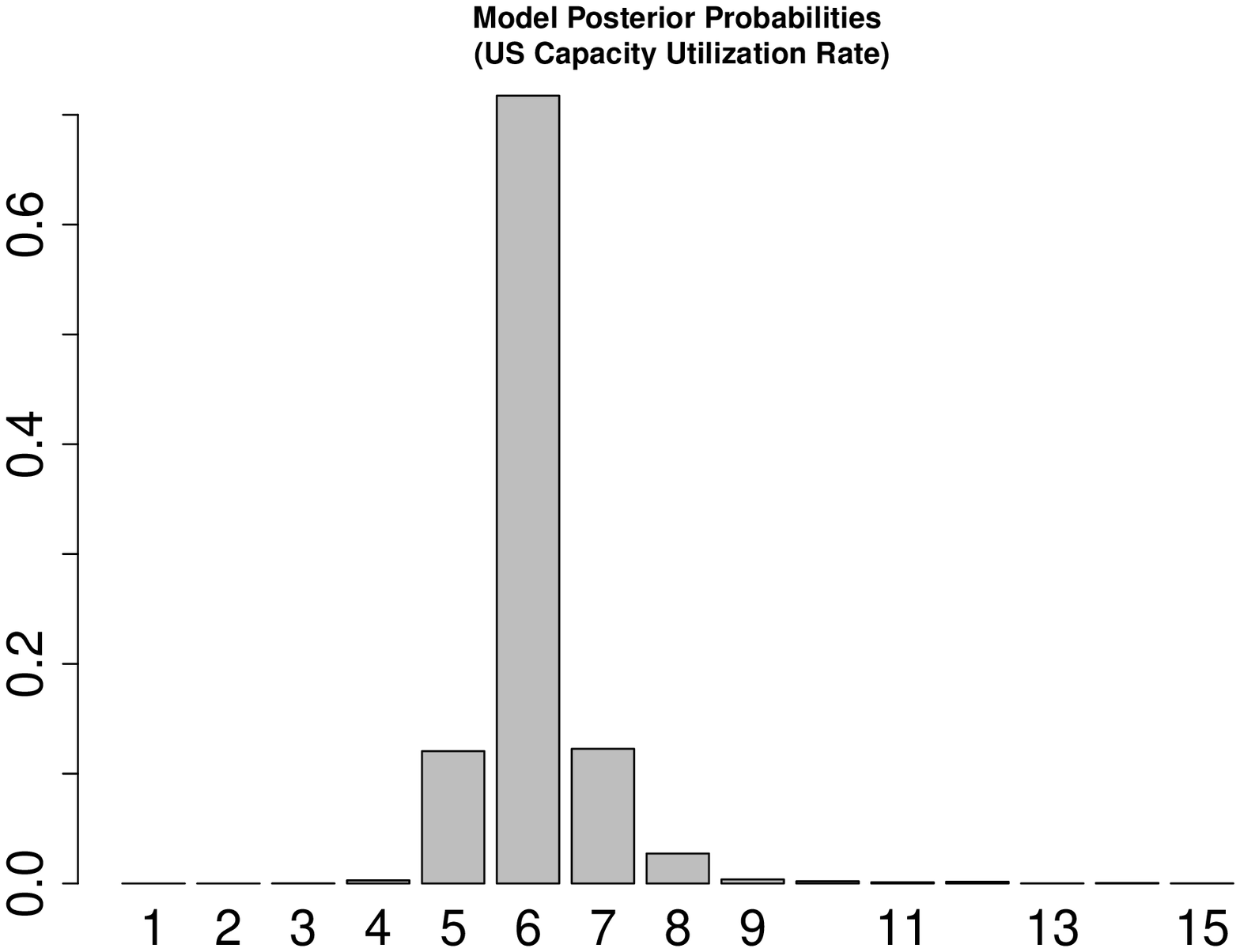}
\caption{Model probabilities for BAR(k) on the US (\textit{upper chart}) and EU13 (\textit{middle chart}) unemployment rates and on the US capacity utilization rate (\textit{bottom chart}).} \label{Fig_BARRJMCMCUnemploy}
\end{figure}

From a graphical inspection (see the bottom chart of Fig. \ref{Fig_BARRJMCMCUnemploy}) we note that the series exhibits a negative trend. A deterministic trend could be naturally included in the Beta regression model with linear conditional mean by imposing some constraints on the slope and the intercept of the linear trend. These constraints can be imposed by a suitable specification of the prior distribution. However in this work we focus on the autoregressive components thus we follow a two steps procedure.

First we define a normalized linear trend $t/T$, where $T$ is the sample size, and introduce the constrained linear regression model
\begin{equation}
x_{t}=\gamma_{0}+\gamma_{1}\frac{t}{T}+\varepsilon_{t},\quad \hbox{with}\,t=1,\ldots,T
\end{equation}
with $\gamma_{0}\in(0,1)$ and $(\gamma_{1}+\gamma_{0})\in(0,1)$. These parameter constraints insure that the residuals of the regression are in the $(0,1)$ interval. In the first step we calculate the de-trended capacity utilization rate $\tilde{x}_{t}=x_{t}-\hat{\gamma}_{1}\frac{t}{T}$. The results of the trend extraction are given in Fig. \ref{Fig_Unempl}.

In the second step we estimate a Beta process on the variable $\tilde{x}_{t}$. The results of the model selection procedure for the Beta process are in Tab. \ref{Tab_Unempl} and Fig. \ref{Fig_BARRJMCMCUnemploy}.

\section{Conclusion}
In this paper we consider Beta autoregressive models of order $k$ for time series data defined on a bounded interval. We focused on conditionally linear mean processes which are particularly suitable for forecasting purposes and allows for an easy interpretation of the parameters in the conditional mean of the process. We proposed a Bayesian procedure in order to estimate the parameters of the model. We showed that, in this context, the Bayesian approach is suitable for dealing with the constraints on the parameters space through a suitable specification of the prior distribution. We introduced three different informative priors and studied, through some simulation experiments, the effects of the prior choice on the inference procedure.

Moreover, we introduced an efficient RJMCMC algorithm for jointly estimating the parameters and selecting the model order. In different sets of simulation
experiments, the proposed algorithm has been shown to be successful in estimating the true parameters for different parameter settings and different choices
of the prior. We also found that the choice of the prior may have effects on the numerical mixing property of the RJMCMC chain by reducing the probability
of the RJMCMC chain moves near the boundaries of the parameter space. In particular, within the three informative priors, the modified truncated Gaussian and
the Beta-type priors could have positive effects on the mixing of the RJMCMC chain.

Through another set of simulation experiments we found that the proposed algorithm is able to find the correct
order $k$ of the autoregressive model. The RJMCMC allowed us to study the effects of the sample size on the model posterior.
In some cases we found evidence of bias and low efficiency of the estimators of the model posterior. Finally we also tested
the performance of the Bayesian inference procedure and of the RJMCMC algorithm by providing two applications to real data.
The first one considers the unemployment rates of United Stases and Euro Area and the second considers the US capacity utilization.

In this paper we focused on the Beta Autoregressive models with conditional linear mean, however the proposed RJMCMC algorithm
can be extended to be applied to other type of Beta processes. In particular the authors are considering the inclusion in the inference process
of the order of the moving average components and the extension to the Beta ARMA processes with nonlinear conditional mean. Both of the extensions
can be considered with and without the inclusion of other explanatory variables.

\newpage

\section*{Appendix A}
When a new reciprocal root $r$ is proposed in the reversible jump MCMC procedure, the relationship
between the autoregressive coefficients of the $\hbox{BAR}(k)$ and of the $\hbox{BAR}(k+1)$ is
\begin{equation}
1-\alpha_{1}^{*}L-\ldots-\alpha_{k+1}^{*}L^{k+1}=\left[1-(1-rL)(1-\alpha_{1}L-\ldots-\alpha_{k}L^{k})\right]
\end{equation}
Identification terms by terms gives the following recursion
\begin{eqnarray}
&&\alpha^{*}_{1}=\alpha_{1}+r\\
&&\alpha^{*}_{j}=\alpha_{j}-r\alpha_{j-1},\quad\hbox{for}\,j=2,\ldots,k+1
\end{eqnarray}
with $\alpha_{k+1}=0$.

The constant term and the autoregressive coefficients belong to the standard simplex, that is
\begin{equation}
\sum_{j=0}^{k}\alpha_{j}\leq 1,\quad\alpha_{l}>0,\, \hbox{for}\, l=0,1,\ldots,k
\end{equation}
for each order $k$ of the Beta process. Thus if the parameters of the
$\hbox{BAR}(k)$ belong to the standard simplex, the following condition
\begin{eqnarray}
&&\alpha^{*}_{0}+\sum_{j=1}^{k+1}\alpha^{*}_{j}\leq1
\Leftrightarrow\alpha^{*}_{0}+(1-r)\sum_{j=1}^{k}\alpha_{j}+r\leq1
\end{eqnarray}
and the positivity constraints for $\alpha_{j}^{*}$, $j=2,\ldots,k$
\begin{equation}
(\alpha_{1}+r)>0,\,(\alpha_{j}-r\alpha_{j-1})>0,\,-r\alpha_{k}\geq0
\end{equation}
are satisfied provided that $$\alpha_{0}^{*}=\alpha_{0}$$ and for $-\alpha_{1}<r\leq0$.

\section*{Appendix B}
For the Gaussian prior we have
\begin{eqnarray}
\nabla^{(1)}\log f(\boldsymbol{\alpha})&=&-\Upsilon^{-1}(\boldsymbol{\alpha}-\boldsymbol{\nu})\\
\nabla^{(2)}\log f(\boldsymbol{\alpha})&=&-\Upsilon^{-1}
\end{eqnarray}
For the modified Gaussian prior we have
\begin{eqnarray}
\nabla^{(1)}\log f(\boldsymbol{\alpha})&=&-\Upsilon^{-1}(\boldsymbol{\alpha}-\boldsymbol{\nu})+
\frac{\kappa}{\phi^{2}\underline{\eta}^{2}(1-\bar{\eta})^{2}}\mathbf{d}\\
\nabla^{(2)}\log f(\boldsymbol{\alpha})&=&-\Upsilon^{-1}+
\frac{\kappa}{\phi^{2}\underline{\eta}^{2}(1-\bar{\eta})^{2}}D-
\frac{2\kappa}{\phi^{2}\underline{\eta}^{3}(1-\bar{\eta})^{3}}\mathbf{d}\mathbf{d}'
\end{eqnarray}
where $D=\boldsymbol{\iota}\mathbf{e}_{1}'+\mathbf{e}_{1}\boldsymbol{\iota}'$, $\mathbf{d}=(\mathbf{e}_{1}-D\boldsymbol{\alpha})'\in\mathbb{R}^{k+1}$, $\boldsymbol{\iota}=(1,\ldots,1)'\in\mathbb{R}^{k+1}$ and $\mathbf{e}_{1}=(1,0,\ldots,0)'\in\mathbb{R}^{k+1}$.

For the transformed Beta prior the $h$-th element of $\nabla^{(1)}\log f(\boldsymbol{\alpha})$ is
\begin{equation}
\partial_{\alpha_{h}}\log f(\boldsymbol{\alpha})=\frac{\nu_{h}-1}{\alpha_{h}}+
\sum_{i=h+1}^{k}\frac{\gamma_{i}+\nu_{i}-1}{A_{i}}-
\sum_{i=h}^{k}\frac{\nu_{i}-1}{A_{i+1}}
\end{equation}
and the $(h,l)$-th element of $\nabla^{(2)}\log f(\boldsymbol{\alpha})$, with $l\geq h$, is
\begin{equation}
\partial_{\alpha_{h}\alpha_{l}}\log f(\boldsymbol{\alpha})=\frac{1-\nu_{h}}{\alpha_{h}^2}-
\sum_{i=l+1}^{k}\frac{\gamma_{i}+\nu_{i}-1}{A_{i}^{2}}+
\sum_{i=l}^{k}\frac{\nu_{i}-1}{A_{i+1}^{2}}
\end{equation}

\bibliographystyle{plainnat}
\bibliography{BetaBib}
\end{document}